 \documentclass[final,authoryear,3p,times,fleqn]{elsarticle}

\usepackage{amssymb}
\usepackage{amsmath}
\usepackage{amsthm}
\usepackage[ruled]{algorithm2e}

\usepackage{graphicx}
\usepackage{float} 

\usepackage{pxfonts,times}
\usepackage{siunitx}

\usepackage{bm}
\usepackage{color}
\usepackage{comment}
\usepackage{natbib}
\usepackage[ruled]{algorithm2e}

\usepackage[
	colorlinks=true, 
	citecolor=blue, 
	urlcolor=blue]{hyperref}

\newtheorem{theo}{Theorem}[section]
\newtheorem{defi}{Definition}[section]
\newtheorem{lemm}{Lemma}[section]
\newtheorem{prop}{Proposition}[section]
\newtheorem{assum}{Assumption}[section]
\newtheorem{coro}{Corollary}[section]

\newtheorem*{prf}{Proof}
\newtheorem*{remark}{Remarks}

\makeatletter
\def\th@plain{\upshape}
\makeatother

\makeatletter
	
	\@addtoreset{equation}{section}
\makeatother

\journal{Elsevier}

\begin{document}

\begin{frontmatter}



\title{Global stability of day-to-day dynamics for\\ schedule-based Markovian transit assignment with boarding queues}



\author[a]{Takashi Akamatsu\corref{cor1}}
\address[a]{Graduate School of Information Sciences, Tohoku University, Miyagi, Japan}
\ead{akamatsu@plan.civil.tohoku.ac.jp}

\author[a]{Koki Satsukawa\corref{cor1}}
\ead{satsukawa@tohoku.ac.jp} 

\author[b]{Yuki Oyama\corref{cor1}}
\address[b]{Department of Civil Engineering, Shibaura Institute of Technology, Tokyo, Japan}
\ead{oyama@shibaura-it.ac.jp}
\cortext[cor1]{Corresponding authors}

\begin{abstract}
Schedule-based transit assignment describes congestion in public transport services by modeling the interactions of passenger behavior in a time-space network built directly on a transit schedule. This study investigates the theoretical properties of scheduled-based Markovian transit assignment with boarding queues. When queues exist at a station, passenger boarding flows are loaded according to the residual vehicle capacity, which depends on the flows of passengers already on board with priority. 
An equilibrium problem is formulated under this nonseparable link cost structure as well as explicit capacity constraints. 
The network generalized extreme value (NGEV) model, a general class of additive random utility models with closed-form expression, is used to describe the path choice behavior of passengers.
A set of formulations for the equilibrium problem is presented, including variational inequality and fixed-point problems, from which the day-to-day dynamics of passenger flows and costs are derived. 
It is shown that Lyapunov functions associated with the dynamics can be obtained and guarantee the desirable solution properties of existence, uniqueness, and global stability of the equilibria. 
In terms of dealing with stochastic equilibrium with explicit capacity constraints and non-separable link cost functions, the present theoretical analysis is a generalization of the existing day-to-day dynamics in the context of general traffic assignment.

\end{abstract}


\end{frontmatter}


\section{Introduction}
In many urban areas, congestion in public transport services remains a major issue, causing boarding queues at stations and on-board discomfort for passengers. 
Congestion in transit networks arises because of the load concentration up to the vehicle capacity in a specific run of a specific line, and passengers make and adjust their choices of route and departure time to minimize their disutility, given their schedule constraints. 
Transit assignment models have been studied to describe such congestion phenomena, and they are either frequency-based or schedule-based. The frequency-based approach represents transit services by lines and computes average on-board flows; thus, it does not explicitly consider single vehicles \citep[e.g.,][]{nguyen1988equilibrium, spiess1989optimal}. 
Currently, advanced transit services operate in a highly regular manner and provide rich information to passengers.
Under such a system, passengers’ route choices should be modeled given the timetable and considering the vehicle capacity of each run \citep{Liu2010}.
The departure time choice behavior also needs to be captured to predict the within-day dynamics of passenger flows. 
The schedule-based approach describes such passengers' choice behavior as a path choice in a time-space network built based on a timetable (hence the approach is also termed \textit{dynamic} transit assignment) and enables us to obtain the loads and performance of each single vehicle \citep[e.g.,][]{Nguyen2001transit, Poon2004transit, Hamdouch2008transit, Sumalee2009transit, Nuzzolo2012, Hamdouch2014}.

From a methodological perspective, four primary challenges arise in dynamic transit assignment: it is necessary to consider for modeling


(1) the heterogeneity of passengers' choice behavior, and 

(2) explicit vehicle capacity constraints. 

\noindent Moreover, such key elements of modeling must be captured within an analytical framework that ensures 

(3) desirable theoretical properties of the solution, and 

(4) computational efficiency. 

\noindent First, to consider the heterogeneous nature of passengers, the probabilistic modeling of path choice behavior is required. 
However, a dynamic transit network is built on a schedule/timetable and contains a complex subnetwork structure to describe interconnections between multiple lines and various transit services. 
Such a network contains a huge number of feasible paths between each origin-destination pair.
Moreover, probabilistic path choice models involve issues of path enumeration and path overlapping. 
Second, to explicitly consider the vehicle capacity of each run, it is required to model the priority of passengers already on-board to boarding passengers and the resultant queues at stations. 
When a queue forms, boarding passengers at the station are loaded according to the residual capacity of vehicle, which is determined based on the flows of precedent on-board passengers. 
In other words, dynamic transit assignment entails an asymmetric link cost structure (i.e., the travel cost on a link is a function of not only the flow on that link but also other link flows). In addition, the discomfort levels of passengers should be differentiated according to their priority.
Third, the analytical formulation of dynamic transit assignment is complex owing to the explicit consideration of capacity constraints. 
Thus, most existing models are simulation-based. 
In simulation-based models, solution convergence is investigated only on a numerical experimental basis, and the theoretical properties of the solution are not clarified. 
From a practical perspective, desirable solution properties such as existence, uniqueness, and stability must be ensured for planners to obtain a reliable solution for equilibrium assignment. 
Finally, to ensure computational efficiency and the desirable theoretical properties, a dynamic transit assignment model should be developed based on an analytical formulation that leads to efficient solution algorithms.

To address the abovementioned four methodological challenges, this study proposes an analytical framework based on Markovian dynamic transit assignment. 
Markovian traffic equilibrium is a link-based stochastic user equilibrium (SUE) model and efficiently incorporates the heterogeneous nature of passengers' choice behavior.
Moreover, its recent extension by \cite{oyama2022markovian} applied in this study can capture path correlation without path enumeration. 
Further, we introduce a novel time-space network representation in which boarding-queue links are explicitly considered to satisfy the residual capacity and first-in-first-out (FIFO) service priority discipline. 
The equilibrium assignment with boarding queues is then mathematically defined. 
To achieve a theoretical analysis, we present a set of formulations for the Markovian dynamic transit assignment, including equivalent variational inequality and fixed-point problems. 
We also develop primal and dual types of day-to-day dynamics that converge to the defined equilibrium under capacity constraints and their Lyapunov functions, according to which we establish the existence, uniqueness, and global asymptotic stability of the equilibrium. 
This day-to-day dynamic addresses a generalized case of \cite{hofbauer2009stable} in terms of having a non-simplex and capacity-constrained solution space. Moreover, the proposed dynamics and Lyapunov functions naturally lead to the development of an efficient solution algorithm, whose convergence are ensured. 

This new framework, for the first time, allows the theoretical analysis of dynamic transit assignment with explicit capacity constraints.
This is a notable contribution to the literature considering that existing models lack the guarantees of desirable solution properties, such as existence, uniqueness and stability.
Moreover, we develop a general class of SUE models, and the investigation of their solution properties by defining day-to-day dynamics and Lyapunov functions also contributes to the general traffic assignment literature.

The remainder of this paper is structured as follows. 
Section 2 reviews the related literature and clarifies the contributions and position of this study.
Section 3 introduces a novel time-space network representation to explicitly consider boarding queues as problem settings. 
Section 4 defines the equilibrium problem.
Section 5 provides a set of formulations, and Section 6 establishes the theoretical properties of the equilibria. 
Finally, Section 7 concludes the paper. 
The appendices provide detailed proofs for the theorems and propositions presented in this study, as well as efficient solution algorithms derived from the theoretical analysis.

\section{Related Literature}

\subsection{Schedule-based transit assignment models}
The transit assignment approach is either frequency-based or schedule-based. 
Although this study focuses on the latter, the reader who is interested in the former approach can refer to a detailed review by \cite{fu2012review}. 

The literature on scheduled-based dynamic transit assignment has often adopted a deterministic user equilibrium (UE) framework instead of SUE~\citep[][]{Ieda1988-ax,Liu2010}. 
In a dynamic transit network that is built based on the timetable and contains a complex subnetwork structure, probabilistic modeling of path choice behavior is a nontrivial task because of path overlapping/correlation and enumeration problems \citep{Hamdouch2014}. 
Logit-based assignment models \citep[e.g.,][]{Nuzzolo2012,cats2016dynamic}, which have been widely applied to SUE, cannot capture the underlying correlation among path utilities because of the independence of the irrelevant alternatives (IIA) property. 
Probit-based models \citep[e.g.,][]{nielsen2006optimisation, Sumalee2009transit} represent the correlation using a covariance matrix but rely on Monte-Carlo simulation that requires extensive computational effort to obtain a good approximation of the expected flow pattern. 
Although path-based models, such as C-logit \citep{Cascetta1996Clogit} and path-size logit/weibit \citep{Xie2020, Wen2021}, capture path overlaps well, they require path enumeration and are inapplicable in a real-size dynamic transit network where the number of feasible paths is almost uncountable. 
Recently, the network generalized extreme value (NGEV) model proposed by \cite{Daly2006NGEV} has been studied in the traffic assignment context \citep{Hara2014NGEV,Ma2015NGEV,oyama2022markovian}, which enables the capture of the path correlation without path enumeration. 
However, such an advanced link-based route choice model has not been applied to dynamic transit assignment.

Dynamic transit assignment models can be further classified into implicit and explicit approaches based on the capacity definition \citep{Nuzzolo2012, Binder2017}.
The implicit approach, similarly to road traffic assignment, defines a non-decreasing cost function with respect to link flow and nominal capacity for each link in the time-space network \citep[e.g.,][]{Tong1999transit, Nuzzolo2001transit}.
Although this simplicity facilitates theoretical analysis of the flow pattern, the implicit approach addresses aggregated link flow, and therefore cannot distinguish boarding passengers at stations from those already on board. 
This approximation in assessing discomfort by ignoring priority is the primary drawback of the implicit approach. 
By contrast, the explicit approach considers the capacity constraint of each vehicle and assigns boarding passengers at stations according to residual capacity. 
Thus, the priority of passengers already on board and the FIFO principle can be captured in this approach \citep[e.g.,][]{Nguyen2001transit, Poon2004transit, Hamdouch2008transit, Sumalee2009transit, Nuzzolo2012, Hamdouch2014}. 
\cite{Binder2017} assigned each passenger to the network in ascending order of an exogenous priority list, and a link becomes unavailable to subsequent passengers when the flow of the link reaches its capacity. \cite{Nuzzolo2001transit,Nuzzolo2012,nuzzolo2016mesoscopic} further investigated the day-to-day dynamics of transit services. 
However, owing to the complexity of the analytical formulation, most existing dynamic transit assignment models with explicit capacity are simulation-based \citep{nuzzolo2009schedule,gentile2016modelling}. 
The development of an analytical model that enables theoretical analysis of solution properties remains challenging. 

\subsection{Theoretical analysis of stochastic user equilibrium}
The theoretical properties of equilibria, such as existence, uniqueness, and stability, under capacity constraints are not very obvious. 
In the literature on traffic assignment, the existence and uniqueness of the solution to a Wardrop deterministic UE problem with capacity constraints have been demonstrated \citep{smith1979existence, smith1983existence, Smith2013}. 
\cite{Smith2013} presented a UE model with explicit link-exit capacities and queuing delays.
Moreover, the existence of a solution was shown under weak conditions. 

The stability of equilibria has been investigated through day-to-day dynamics that capture the evolution of network flow patterns over time. 
In particular, continuous-time dynamical systems that converge to deterministic UE have been extensively studied, some of which established general conditions that ensure global asymptotic stability of the solution \citep[e.g.,][]{Smith1984, friesz1994day, zhang1996local, Han2012, Guo2015}.
\cite{Yang2009} investigated a general dynamic called rational behavior adjustment process that includes many deterministic dynamics as special cases. 

While deterministic models assume that travelers possess complete information about the network, the dynamical systems that converge to the SUE assume that travelers make heterogeneous perceptual errors, and thus are more realistic.
The literature on the dynamics for SUE can be divided into discrete-time (perception-based) and continuous-time (flow-based) models.
The relationships between stability/instability properties in discrete- and continuous-time models were explored by \cite{Watling1999} and \cite{Cantarella2016a}. 
For mathematical convenience, the theoretical properties of SUE have been mainly investigated within continuous-time day-to-day dynamics, including the logit dynamic \citep{Watling1999, sandholm2010population} and logit-based Smith dynamic \citep{Smith2016a} whose convergence and stability are established.

\cite{Xiao2019a} recently proposed a general framework for continuous day-to-day dynamics converging to SUE, and the logit dynamic and logit-based Smith dynamic were contained as special cases in the framework. 
They used \cite{fisk1980some}'s equivalent optimization formulation of SUE as a general Lyapunov function, which guarantees the Lyapunov stability of the framework. 
In their study, separable link cost functions are assumed because Fisk's equivalent optimization formulation is not available for an asymmetric case in which cost functions are nonseparable.
\cite{Guo2013a} established the convergence for a discrete-time form of the link-based logit dynamic in the appendix of their paper. 
Under separable link cost conditions, they used the objective function of the equivalent optimization of the unconstrained program of \cite{Sheffi1985} for the Lyapunov function and analyzed the SUE. 

The logit-based Smith dynamic proposed by \cite{Smith2016a} is an extension of the logit dynamic based on the logic of the Smith dynamic \citep{Smith1984}, and its global stability has been established under more generalized link cost conditions.
Further, \cite{ye2022stochastic} presented day-to-day SUE models whose stability is ensured under non-separable link cost conditions. 
The author transformed a logit-based perceived cost dynamic in continuous-time into a flow-based dynamic and showed its equivalence to the logit-based extension of \cite{Jin2007}'s dynamic. 
By Using this equivalent path flow-based formulation, \cite{ye2022stochastic} proved the global stability of the logit-based perceived cost dynamic. 
The author additionally introduced other types of route choice models, such as C-logit, path-size logit, and Weibit models, and examined different equilibria and global stability of different dynamics with a non-separable link cost function. 
However, existing dynamics with global stability ensured under non-separable link cost conditions are all path-based, and thus suffer from the problems of initial path flow patterns and path overlapping \citep{He2010}.

\subsection{Contributions of the study}
Given the literature review, we state that existing dynamic transit assignment models with explicit capacity constraints are mostly based on a deterministic UE framework, which ignores the heterogeneous/stochastic nature of the path choice behavior of passengers. 
Theese are also simulation-based models whose solution properties, such as existence, uniqueness, and stability, have not been elucidated. 
An analytical framework that enables rigorous theoretical analysis of the solution properties remains undeveloped for dynamic transit assignment. 
Moreover, in the literature on day-to-day dynamical systems, the global stability has been established for either deterministic or under separable link cost conditions. 
\cite{ye2022stochastic} is a notable exception in which the global stability is established for dynamical systems for SUE with non-separable link cost functions.
However, the  Lyapunov function is based on flows at equilibrium and cannot be used to develop a solution algorithm. 
\cite{Smith2016a} deals with generalized link cost conditions, but their model is logit-based and cannot capture the path overlapping/correlation.
Moreover, no studies have investigated the global stability of dynamical systems for SUE under strict capacity constraints.

This study develops an analytical framework for dynamic transit assignment. 
The contributions lie in both modeling and theoretical analysis, which are summarized as follows:
\begin{itemize}
    \item \textbf{Markovian dynamic transit assignment model}. 
    This study proposes a link-based SUE model for dynamic transit assignment with explicit consideration of capacity constraints and boarding queues. 
    We propose a novel time-space network representation in which boarding-queue links are explicitly introduced to satisfy the residual capacity and FIFO service priority discipline. 
    In this dynamic transit network, we formulate the Markovian traffic equilibrium based on \cite{oyama2022markovian}.
    This model is a state-of-the-art and link-based SUE model that generalizes closed-form GEV models and captures the correlation among time-space paths without path enumeration.
    Thus, the heterogenous/stochastic nature of the passengers' choice behavior is captured while retaining computational efficiency. 
    
    \item \textbf{Theoretical analysis of SUE under capacity constraints}. 
    To achieve the theoretical analysis, we present a set of formulations, including equivalent variational inequality and fixed-point problems. 
    Two types of generalized day-to-day dynamics for SUE are derived from the formulations: one is based on link flows and the other on link costs. 
    Lyapunov functions associated with the dynamics are developed, and the existence, uniqueness, and global stability of the equilibrium are established. 
    The dynamic and Lyapunov function defined based on link flows deal with a generalized case of \cite{hofbauer2009stable} (Theorem 5.2), which has never been applied or extended in the traffic assignment context in terms of having a non-simplex and capacity-constrained solution space.
    Moreover, the proposed Lyapunov function has several advantages: (a) it is compatible with any additive random utility (RU) model; (b) it naturally leads to the development of an efficient solution algorithm whose convergence is ensured; (c) it is link-based and does not suffer from the problems of initial path flow pattern and path overlapping; and (d) it can consider explicit capacity constraints.
    Thus, this contribution on theoretical analysis is not limited to dynamic transit assignment but is also beneficial to the general traffic assignment context.
\end{itemize}

\section{Problem settings}\label{Sec:Settings}

\subsection{Time-space network representation with boarding queues}
Consider a fixed timetable of public transport and fixed passenger demand. 
We develop a dynamic transit assignment model to predict passenger flow in sections of train-runs and queues at stations. 
To this end, the path choice behavior of passengers is described in a directed time-space network $G = (\mathcal{N}, \mathcal{L})$, a graph representation of the timetable, where $\mathcal{N}$ and $\mathcal{L}$ are the sets of nodes and links.
First, we construct a time-space network as in the literature \citep[e.g.,][]{Binder2017}: Let $\mathcal{T}$, $\mathcal{S}$, and $\mathcal{R}$ be the sets of discrete time steps, transit stations, and train-runs, respectively. 
A time-space node $i = (t, s, r) \in \mathcal{N}$ is characterized by a combination of time $t \in \mathcal{T}$, station $s \in \mathcal{S}$, and run $r \in \mathcal{R}$. 
A link $(i,j) \in \mathcal{L}$ connects two nodes $i, j \in \mathcal{N}$ and represents the physical/virtual time-space movement of passengers.
Additionally, we consider dummy nodes for origins and destinations.
An origin dummy node $o = s_o \in \mathcal{O} \subset \mathcal{N}$ represents a departure station $s_o$, and a destination dummy node $d = (t_d, s_d) \in \mathcal{D} \subset \mathcal{N}$ represents an arrival station $s_d$ with arrival time $t_d$. 
The given passenger demand $\mathbf{q} \equiv \{q^d_o\}_{o \in \mathcal{O}, d \in \mathcal{D}}$ is differentiated by origin-destination (OD) pair, that is, the departure and arrival stations, as well as the desired arrival time. 
The set of OD pairs is defined as ${\mathcal W} \equiv \mathcal{O} \times \mathcal{D}$.

\begin{figure}[t]
	\centering
	\includegraphics[width=0.8\linewidth]{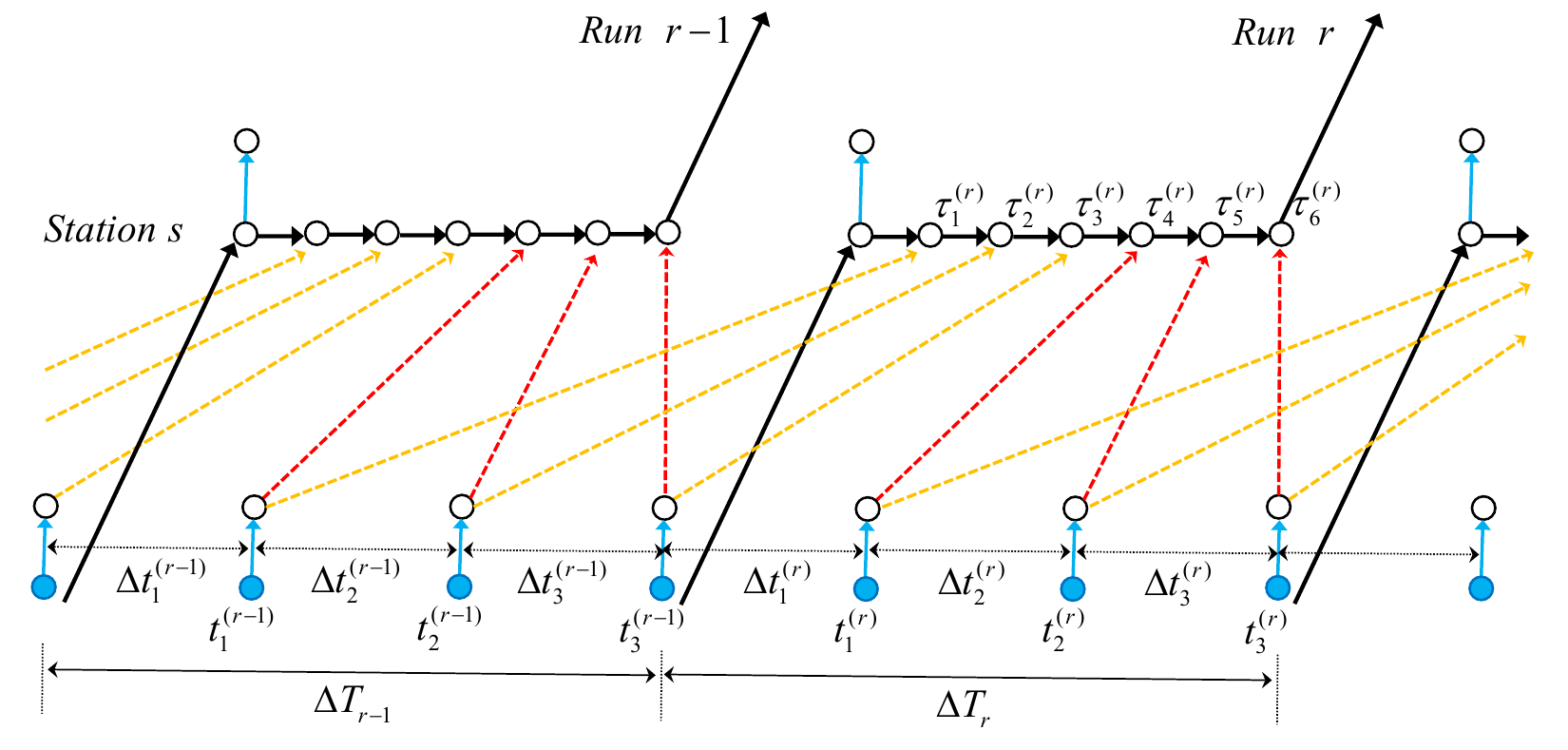}
	\caption{Time-space sub-network for station $s$; 
	blue nodes represent source (origin) nodes of station $s$, where $K=3$ is assumed. 
	The red arrows entering into each run $r$ denote ``boarding-queue links," where $M=2$ is assumed.
	For notational brevity, the subscript $s$ representing a station is omitted.}
	\label{Fig:TSNetwork}
\end{figure}

In this study, to explicitly consider boarding queues at stations, we extend the time-space network by introducing a subnetwork for each station-run pair $(s, r)$, $s \in S, r \in R$ (Figure~\ref{Fig:TSNetwork}).
In the subnetwork representation, two new types of nodes are introduced: ``source'' and ``on-board'' nodes. 
A source node represents an event of passenger arrival at the station from the origin node or by transfer. 
Let $T_{s,r}$ be the departure time of train-run $r$ from station $s$. 
We divide the headway $\Delta T_{s,r} \equiv [T_{s,r-1}, T_{s,r}]$ of two consecutive runs into $K$ periods. 
Passengers arriving at station $s$ during the $k$-th period $\Delta t^{(s,r)}_{k}$ are represented as flows emanating from a source node of station $s$ at time $t^{(s,r)}_k$. 
In addition, the source nodes are connected by access links from the origin for station $s$, describing the departure time choice behavior of passengers.

An on-board node represents the on-board state of passengers at station $s$ on train-run $r$. 
We assume that $M$ boarding queues can be formed simultaneously at each station, where the $m$-th boarding queue at a given point in time is served by the $m$-th train that arrives after that time. 
Passengers arriving at each station choose one of the queues to wait for a run, where the service discipline for each queue is FIFO.
In the subnetwork representation, a single queue waiting for run $r$ at station $s$ is divided into $L \equiv K \times M$ ``boarding-queue'' links (see Figure~\ref{Fig:TSNetwork}), each of which emanates from one of source nodes for the station $s$ and enters one of on-board nodes. 
For each station-run pair $(s,r)$, the $l$-th on-board node has an associated boarding time $\tau^{(s,r)}_l$ of the passenger flow on the $l$-th boarding-queue link, and the on-board nodes are arranged in ascending order with respect to the boarding times $\{\tau^{(s,r)}_l\}$. 
To satisfy the FIFO service discipline, the boarding-queue links for each run do not cross each other.
The arrival time at the source node of passengers who use the $l$-th boarding-queue link must be earlier than that for the $l+1$-th boarding queue link. 
As explained later, this representation also explicitly describes the priority of passengers already on-board on run $r$ to passengers boarding run $r$ from station $s$.

The set of links $\mathcal{L} = \mathcal{L}_{T} \cup \mathcal{L}_{A} \cup \mathcal{L}_{B} \cup \mathcal{L}_Q$ is classified into the subsets of time-space (driving) links $\mathcal{L}_T$, access/egress/transfer links $\mathcal{L}_A$, on-board links $\mathcal{L}_B$, and boarding-queue links $\mathcal{L}_Q$. Each of these subsets is further divided into mutually exclusive subsets using station-run pairs; for instance,
\begin{align}
	\mathcal{L}_{Q} = \cup_{(r,s)\in \Psi}\mathcal{L}_{Q}^{(r,s)}
	\ \text{and}\ 
	\mathcal{L}_{Q}^{(r,s)}\cap \mathcal{L}_{Q}^{(r',s')} = \emptyset,
	\quad \forall (r,s)\neq (r',s')\in \Psi,
\end{align}
where $\Psi$ denotes the set of station-run pairs.

\subsection{Travel cost and boarding-queues with priority}

\begin{assum}\label{Assum:Monotonicity}
	The travel cost $c_{ij}$ of each link $(i,j)\in\mathcal{L}_{A}\cup \mathcal{L}_{B}$ is given by a continuous, differentiable and strictly increasing function of the link flow $x_{ij}$.
\end{assum}

\begin{figure}[t]
	\centering
	\includegraphics[width=0.6\linewidth]{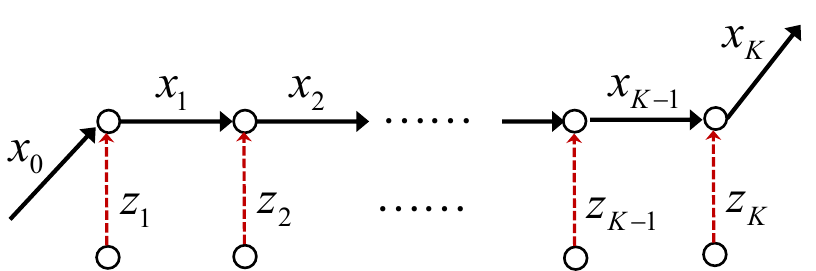}
	\caption{Flows on ``on-board" and ``boarding-queue" links for each run-station pair; 
	the black/red arrows denote on-board/boarding-queue links. 
	(For simplicity, the superscript $(r,s)$ is omitted.)}
	\label{Fig:FlowDescription}
\end{figure}

For clarity in presenting the travel cost for links in $\mathcal{L}_{Q}$ (boarding-queue links), we use local notation here.
For each run-station pair $(r,s)\in\Psi$, let $\mathbf{x}^{(r,s)}$ and $\mathbf{z}^{(r,s)}$ be the flows on ``on-board links" and ``boarding-queue links," respectively.
As can be seen from the structure of the time-space network (see Figure~\ref{Fig:FlowDescription}), the link flows $\mathbf{x}^{(r,s)}$ and $\mathbf{z}^{(r,s)}$ satisfy the following flow conservation:
\begin{align}
	x_{l}^{(r,s)} 
	&= x_{l-1}^{(r,s)} + z_{l}^{(r,s)}\\
	&= x_{0}^{(r,s)} + \sum_{k=1}^{l}z_{k}^{(r,s)},\quad (l=1,2,\ldots,L).\label{Eq:FlowConservBasis}
\end{align}

\begin{assum}\label{Assum:InftyCost}
	Each on-board link for run $r$ has a capacity $\mu^{(r)}$.
	The travel cost of the $k$-th boarding-queue link for the run-station pair $(r,s)$ is an increasing function of the flow $x_{l}^{(r,s)}$ of the successive on-board link:
	\begin{align}
		C_{Q(l)}^{(r,s)}\equiv f_{l}(x_{l}^{(r,s)}\mid \mu^{(r)}),\label{Eq:OnBoardCost}
	\end{align}
	where $f_{l}(x\mid \mu)$ is a continuous, differentiable, and strictly increasing function of $x$ that satisfies 
	\begin{align}
		\lim_{x\to \mu}.f_{l}(x\mid \mu) = +\infty.\label{Eq:InftyCost}
	\end{align}
\end{assum}

\begin{assum}\label{Assum:SelectFiniteCost}
	Each passenger always chooses a boarding-queue link with a finite cost.
\end{assum}

Eq.~\eqref{Eq:InftyCost} in \textbf{Assumption~\ref{Assum:InftyCost}} together with \textbf{Assumption~\ref{Assum:SelectFiniteCost}} implies that the flow on the $k$-th on-board link should satisfy $x_{l}(r,s)<\mu^{(r)}$.
It follows from this and the flow conservation~\eqref{Eq:FlowConservBasis} that the flow $z_{l}(r,s)$ on the $k$-th boarding-queue link should be adjusted (by users' rational choice behaviour) so as to satisfy $x_{l-1}^{(r,s)}+z_{l}^{(r,s)}<\mu^{(r)}$, or equivalently,
\begin{align}
	z_{l}^{(r,s)}<\mu^{(r)} - x_{l-1}^{(r,s)}.\label{Eq:FlowLimit}
\end{align}
The right-hand side of \eqref{Eq:FlowLimit}, $\mu^{(r)} - x_{l-1}^{(r,s)}$ means the ``residual capacity" of the train-run $r$ left after passengers in the $(l-1)$-th boarding-queue have boarded.
In other words, the inequality~\eqref{Eq:FlowLimit} implies that the flows on the two consecutive boarding-queue links are consistent with the FIFO discipline.
Because the same logic holds for all $l=1,2,\ldots,L$ and all run-station pairs, under \textbf{Assumptions~\ref{Assum:InftyCost}} and \textbf{\ref{Assum:SelectFiniteCost}}, all boarding-queue link flows are consistent with the FIFO discipline and all on-board link flows are always within the capacity limitations.

Substituting the flow conservation~\eqref{Eq:FlowConservBasis} into \eqref{Eq:OnBoardCost}, we obtain the boarding-queue cost represented as a function of $\mathbf{z}^{(r,s)}$ and $x_{0}^{(r,s)}$:
\begin{align}
	C_{Q(l)}^{(r,s)}\equiv f_{l}(x_{0}^{(r,s)} + \sum_{k=1}^{l}z_{k}^{}(r,s) \mid \mu^{(r)}),
\end{align}
Note that the cost of the $k$-th boarding-queue link is a function of the flows of boarding-queue links upstream of the $l$-th link (i.e., a function of flows with earlier arrival times).
This implies that the link cost vector $\mathbf{c}_{Q}^{(r,s)}$ for each run-station pair $(r,s)$ is a function of $\mathbf{z}^{(r,s)}$ with a lower triangular Jacobian.
This observation leads us to the following lemma:
\begin{lemm}\label{Lemm:Monotone_BoardingLink}
	Suppose that $x_{0}^{(r,s)}+\sum_{k=1}^{l}z_{k}^{(r,s)}<\mu^{(r)}$ $(l=1,2,\ldots,L)$.
	Then, the vector cost function $\mathbf{c}_{Q}^{(r,s)}$ of the links in $\mathcal{L}_{Q}^{(r,s)}$ for each run-station pair $(r,s)$ is strictly monotone with respect to $\mathbf{z}^{(r,s)}$.
\end{lemm}
\begin{prf}
	For brevity, we omit the superscript $(r,s)$ for denoting the run-station pair.
	Note that the costs of boarding-queue for the run-station pair $(r,s)$ depend only on the flows $\mathbf{z}$ of the same run-station pair.
	The Jacobian $\mathbf{J}_{Q}$ of the cost function $\mathbf{C}_{Q}$ is given by the following lower triangular matrix:
	\begin{align}
		[\mathbf{J}_{Q}]_{i,j}\equiv\partial C_{Q(i)}/\partial z_{j}
		=
		\begin{cases}
			f'_{i}(x_{i})\quad &\text{if}\ i\geq j\\
			0					&\text{if}\ i< j
		\end{cases}
	\end{align}
	Herein, all the diagonal elements (eigenvalues) of $\mathbf{J}_{Q}$ are always positive according to \textbf{Assumption~\ref{Assum:InftyCost}}.
	This implies that $\mathbf{J}_{Q}$ is positive definite; therefore, $\mathbf{C}_{Q}(\mathbf{z})$ is strictly monotone.\qed
\end{prf}

We then define the following supply-feasibility, ensuring that a flow pattern is within the capacity: 
\begin{defi}
	The set $Y$ of the supply-feasible flows of the time-space network is defined as the link flows $\mathbf{x}\in\mathcal{R}_{+}^{L}$ for which $\mathbf{C}(\mathbf{x})$ is finite.
	Set $C$ of feasible link cost patterns is defined as 
	\begin{align}
		C\equiv\{\mathbf{c}\in\mathcal{R}_{+}^{L}\mid \mathbf{c}=\mathbf{C}(\mathbf{x}),\ \forall \mathbf{x}\in Y \}.
	\end{align}
\end{defi}

\textbf{Assumption~\ref{Assum:Monotonicity}} suggests that the partial derivative of the cost function of a link $l\in \mathcal{L}_{A}\cup \mathcal{L}_{B}$ with respect to a link flow, $\partial c_{l}(x_{l})/\partial x_{l'}$, is positive if $l=l'$ or zero otherwise, for a given flow pattern within $Y$.
This assumption and \textbf{Lemma~\ref{Lemm:Monotone_BoardingLink}} indicate that the link cost vector $\mathbf{C}(\mathbf{x})$ is a function of $\mathbf{x}$ with a lower triangular Jacobian whose diagonal elements are always positive.
Hence, we derive the following propositions:
\begin{prop}\label{Prop:CostMonotone}
	The vector link cost function $\mathbf{C}:Y\rightarrow C$ of the entire time-space network for the NGEV transit assignment is strictly monotone with respect to link flows $\mathbf{x}\in Y$; that is
	\begin{align}
		(\mathbf{C}(\mathbf{x}') - \mathbf{C}(\mathbf{x}''))\cdot (\mathbf{x}'-\mathbf{x}'')>0,\quad \forall \mathbf{x}'\neq \mathbf{x}''\in Y.
	\end{align}
\end{prop}

\begin{prop}\label{Prop:InverseCost}
	The vector link cost function $\mathbf{C}$ is a one-to-one function of $\mathbf{x}\in Y$ and has an inverse function $\mathbf{C}^{-1}:C\rightarrow Y$, which is strictly monotone with respect to link costs $\mathbf{c}\in C$; that is,
	\begin{align}
		(\mathbf{C}^{-1}(\mathbf{c}') - \mathbf{C}^{-1}(\mathbf{c}''))\cdot (\mathbf{c}'-\mathbf{c}'')>0,\quad \forall \mathbf{c}'\neq \mathbf{c}''\in C.
	\end{align}
\end{prop}

\section{Network GEV equilibrium assignment}\label{Sec:NetworkGEV}
This section formulates a network generalized extreme value (NGEV) equilibrium assignment problem that includes the dynamic transit assignment problem defined in the previous section.
We first introduce the formulation of the NGEV route choice model.
we next formulate a \textit{flow-independent} NGEV assignment model, which is a network loading model based on the NGEV route choice model.
We then present a \textit{flow-dependent} NGEV equilibrium assignment problem based on the flow-independent NGEV assignment model.

\subsection{NGEV route choice model}
The NGEV model, which was proposed by \cite{Bierlaire2002-tn,Daly2006NGEV}, is one of the most flexible closed-form models for representing utility correlation across alternatives by constructing a GEV network.
In the context of traffic assignment, \cite{Papola2013-ol}, \cite{Hara2014NGEV}, and \cite{Mai2016-zr} proposed route choice models based on the NGEV model.
This NGEV route choice model directly utilizes the transportation network structure as a GEV network to capture the underlying correlation among route utilities, and provides route choice probabilities compliant with the framework of RU models without route enumeration\footnote{For details of the NGEV model, see for example, \cite{Daly2006NGEV} and \cite{Papola2013-ol}.}.
\cite{oyama2022markovian} recently established a framework of the traffic equilibrium assignment based on the NGEV route choice model.

Given an OD pair $(o,d) \in {\mathcal W}$, a route choice can be described as a joint choice of all the elemental links of the route in their ordered sequence from origin $o$ toward destination $d$.
In this \textit{ordered joint choice context}~\citep[][]{Papola2013-ol}, the choice probability $p(r)$ of a route $r$ is given by the product of the conditional choice probabilities of all the links on the route as follows:
\begin{align}
	p(r) = \prod_{ij \in \mathcal{L}_r} p^d_{ij|i},
	\label{Eq:RouteProb}
\end{align}
where $p^d_{ij|i}$ is the conditional choice probability of link $ij$ on node $i$ for a traveler toward destination $d$; $\mathcal{L}_{r}$ is the set of links on route $r$.
The conditional choice probability in the NGEV model is defined as follows~\citep[][]{Daly2006NGEV}:
\begin{align}
	p^d_{ij|i} = \frac{\alpha^d_{ji} (G^d_{ij})^{\theta^d_i/\theta^d_j}}
	{{\displaystyle \sum_{j' \in \mathcal{F}(i)}} \alpha^d_{j'i} (G^d_{ij'})^{\theta^d_i/\theta^d_{j'}}},
	\label{Eq:PijG}
\end{align}
where $\mathcal{F}(i)$ is the set of successor nodes for node $i$; 
$\theta^d_{i}$ is the variance scale parameter associated with node $i\in\mathcal{N}$, which is strictly positive;
$\alpha^d_{ji}$ is the allocation parameter that represents the degree of membership of node $j\in\mathcal{N}$ to the predecessor node $i \in \mathcal{B}(j)$, satisfying $\sum_{i \in \mathcal{B}(j)} \alpha^d_{ji} = 1$, $\alpha^d_{ji} \ge 0$.
$G^d_{ij}$ is the NGEV generating function associated with link $ij\in\mathcal{L}$.
This function can be related to the expected minimum cost $s_{id}$ from node $i$ to destination $d$ in the following equation~\citep[see,][]{oyama2022markovian}:
\begin{equation}
	G^d_{ij} =  e^{- \theta^d_j (c_{ij} + s_{id} )}.
	\label{Eq:Gmu}
\end{equation}

By combining Eqs.~\eqref{Eq:Gmu} and \eqref{Eq:PijG}, the NGEV route choice model for a given $\mathbf{c}\in C$ is given as :
\begin{align}
	&p_{ij\mid i}^{d} = 
	\cfrac{\alpha_{ji}^{d}e^{-\theta_{i}^{d}(c_{ij} + s_{jd})}}{\sum_{n\in \mathcal{F}(i)}\alpha_{ni}^{d}e^{-\theta_{i}^{d}(c_{in} + s_{nd})}}
	&&\forall ij\in\mathcal{L},\forall d\in \mathcal{D},
		\label{Eq:LinkChoiceProb}\\
	&1 = \sum_{n\in\mathcal{F}(i)}\alpha_{ni}^{d}e^{-\theta_{i}^{d}(c_{in} + s_{nd} - s_{id})}
	&&\forall i\in\mathcal{N},\forall d\in\mathcal{D}.
		\label{Eq:LinkChoiceProbSum}
\end{align}

%
%

\subsection{Formulation of NGEV assignment}
Let $\mathbf{q}_{d}\equiv (q_{i}^{d})_{i\in\mathcal{N}}$ be a given demand flow toward a destination $d\in\mathcal{D}$ from every node in $\mathcal{N}$.
We assume that $q_{i}^{d} > 0$ for all $(i,d)\in\mathcal{W}$ and $q_{i}^{d} = 0$ otherwise.
For the given demand flow, the NGEV assignment derives link flow based on the NGEV route choice model introduced in the previous section.
It consists of link flow distribution and flow conservation conditions.

The former condition describes the link flow pattern resulting from the NGEV route choice model.
Let $x_{ij}^{d}\geq 0$ be the non-negative disaggregated link flow toward a destination $d\in\mathcal{D}$ on a link $ij\in\mathcal{L}$.
Then, this link flow satisfies the following condition:
\begin{align}
	x_{ij}^{d} = p_{ij\mid i}^{d}\left( q_{i}^{d} + \sum_{n\in \mathcal{B}(i)}x_{ni}^{d}\right),\label{Eq:LFDist_Proto}
\end{align}
where $p_{ij\mid i}$ is a link choice probability defined in Eq.~\eqref{Eq:LinkChoiceProb}.

Next, for a disaggregated link flow toward a destination $d\in\mathcal{D}$, the following flow conservation condition should be satisfied at each node $i\in\mathcal{N}\setminus \{d\}$\footnote{Since the flow conservation condition at the destination $d$ can be derived from the conditions at the other nodes (i.e. redundant), we omit it.}:
\begin{align}
	\sum_{n\in \mathcal{F}(i)}x_{in}^{d} - \sum_{n\in \mathcal{B}(i)}x_{ni}^{d} = q_{i}^{d}.\label{Eq:FCNodeElement}
\end{align}
This condition can be expressed in the following vector form for every $d\in\mathcal{D}$:
\begin{align}
	\mathbf{A}_{d}\mathbf{x}_{d} = \mathbf{q}_{d},\label{Eq:FCNode}
\end{align}
where $\mathbf{A}_{d}$ denotes a reduced incident matrix $(|\mathcal{N}-1|\times |\mathcal{L}|)$.
This matrix is constructed by eliminating the row correspond to the destination $d$ of a node-link incident matrix $\mathbf{A}\equiv [a_{k,ij}]_{k\in \mathcal{N}, ij\in\mathcal{L}}$: $a_{k,ij} = 1$ if link $ij$ is leaving from node $k$, $a_{k,ij} = -1$ if link $ij$ is arriving at node $k$ and $a_{k,ij} = 0$ otherwise.
$\mathbf{x}_{d}\equiv (x_{ij}^{d})_{ij\in\mathcal{L}}$ is the vector of the link flow toward the destination $d\in\mathcal{D}$.

By substituting Eq.~\eqref{Eq:FCNodeElement} into Eq.~\eqref{Eq:LFDist_Proto}, the link flow distribution condition can be expressed as follows:
\begin{align}
	&x_{ij}^{d} = p_{ij\mid i}^{d}\sum_{n\in\mathcal{F}(i)}x_{in}^{d},
	&&\forall ij\in\mathcal{L},\forall d\in\mathcal{D}.\label{Eq:LFDist}
\end{align}
This equation is a simple representation of the forward type link choice conditional probability from the disaggregated link flow.

Here, we define a vector variable $\mathbf{x}\equiv (\mathbf{x}_{d})_{d\in\mathcal{D}}$ and a set $\mathbf{s}\equiv\{\mathbf{s}_{d},\forall d\in\mathcal{D}\}$, where $\mathbf{s}_{d}\equiv(s_{id})_{i\in\mathcal{N}}$.
The (flow-independent) NGEV assignment is then defined as follows:
\begin{defi}[\textit{NGEV assignment}]
	For a given link cost pattern $\mathbf{c}\in C$, the NGEV assignment is mathematically defined as the problem of finding a solution tuple $\langle \mathbf{x}, \mathbf{s}\rangle$ that satisfies Eqs.~\eqref{Eq:LinkChoiceProb}-\eqref{Eq:LFDist}.
\end{defi}

Note that $s_{id}$ is derived as a function of the given link cost pattern $\mathbf{c}$ through the following recursive formulation~\citep[][]{oyama2022markovian}:
\begin{align}
&
	s_{id}\equiv S_{id}(\mathbf{c}) 
	=  
	-\frac{1}{\theta^d_i} \ln \sum_{j \in \mathcal{F}(i)} \alpha^d_{ji} e^{-\theta^d_i (c_{ij} + S_{jd}(\mathbf{c}))},
	&&
	\forall i\in\mathcal{N},\forall d\in\mathcal{D}.
	\label{Eq:NGEV_ExpMinCost}
\end{align}

\subsection{NGEV equilibrium assignment}
Based on the NGEV assignment, we now define the flow-dependent NGEV equilibrium assignment.
First, we define the relationship between the disaggregated and aggregated link flow as follows:
\begin{align}
	\mathbf{x} = \sum_{d\in\mathcal{D}}\mathbf{x}_{d}.
\end{align}
We further define the following set of \textit{demand-feasible aggregated} flow:
\begin{defi}[\textit{Demand-feasible flow}]
	The set $X$ of demand-feasible aggregated link flows $\mathbf{x}$ is defined as 
	\begin{align}
	&X\equiv \{ \mathbf{x} \in \mathbb{R}_{+}^{L} \mid \mathbf{x} = \sum_{d\in\mathcal{D}}\mathbf{x}_{d},\ \mathbf{x}_{d}\in X_{d}, \forall d\in\mathcal{D} \},\label{Eq:DemandFeasible}
	\end{align}
	where $X_{d}$ is the set of demand-feasible disaggregated flows satisfying Eq.~\eqref{Eq:FCNode}.
	
	The set $\hat{X}$ of demand-feasible disaggregated link flows $\hat{\mathbf{x}}$ is also defined as
	\begin{align}
		\hat{X}\equiv \prod_{d\in\mathcal{D}}X_{d}=
		\{ \hat{\mathbf{x}}=[\mathbf{x}_{1},\ldots,\mathbf{x}_{D}]\in \mathbb{R}_{+}^{L\times D} \mid \mathbf{x}_{d}\in X_{d}\ \forall d\in\mathcal{D} \}.\label{Eq:DemandFeasibleDisaggregate}
	\end{align}
\end{defi}

Then, the NGEV equilibrium assignment is mathematically defined as follows:
\begin{defi}[\textit{Network GEV equilibrium assignment}]
The NGEV equilibrium assignment is mathematically defined as the problem of finding a solution tuple $\langle \mathbf{c}, \mathbf{x},\mathbf{s} \rangle$ that satisfies Eqs.~\eqref{Eq:LinkChoiceProb}-\eqref{Eq:DemandFeasibleDisaggregate} and the conditions of the link cost function defined in Section~\ref{Sec:Settings}.
\end{defi}
\noindent It should be noted that our problem differs from that of \cite{oyama2022markovian} in that the capacity constraints (i.e., the supply feasibility of link flows) are explicitly considered.

\section{Equivalent fixed-point and variational inequality problems}
This section presents several equivalent problems to the NGEV equilibrium assignment that arise naturally from this problem.
First, we introduce a conjugate function of the generalized entropy function as a basis for the analysis of the theoretical properties of the equilibrium assignment problem.
We then construct equivalent fixed-point and variational problems using this function.

\subsection{Conjugate functions}
First, we introduce the following lemma showing an equivalent optimization problem:
\begin{lemm}[\textbf{Proposition 1} and \textbf{Corollary 1} in \cite{oyama2022markovian}]
For a given link cost pattern $\mathbf{c}\in C$, the flow-independent NGEV assignment for a network with a many-to-one OD demand pattern is equivalent to the following convex programming having a unique solution:
\begin{align}
	&	\min_{\mathbf{x}_{d}\in X_{d}}Z_{d}(\mathbf{x}_{d})\equiv 
	\mathbf{c}\cdot \mathbf{x}_{d} - H_{d}(\mathbf{x}_{d}),\label{Eq:CP-NGEV1}
	\\
	&	\text{where}\quad 
	H_{d}(\mathbf{x}_{d})\equiv -\sum_{i\in\mathcal{N}}\cfrac{1}{\theta_{i}^{d}}\sum_{j\in\mathcal{F}(i)}x_{ij}^{d}\left[ \ln (x_{ij}^{d}/\sum_{j\in\mathcal{F}(i)}x_{ij}^{d}) - \ln \alpha_{ij}^{d} \right].\label{Eq:CP-NGEV2}
\end{align}
\end{lemm}
\noindent An entropy function $H_{d}$ for NGEV assignment for a destination $d\in\mathcal{D}$ is defined by a disaggregated link flow pattern $\mathbf{x}_{d}$.
It is worth noting that this generalized entropy function is bounded because each disaggregated link flow $x_{ij}^{d}$ is bounded by the flow conservation condition~\eqref{Eq:FCNode}, and the argument $x_{ij}^{d}/\sum_{j\in \mathcal{F}(i)}x_{ij}^{d}$ is bounded with $0$ and $1$ by definition.

In association with this problem, we define the optimal solution map $\mathbf{F}_{d}(\mathbf{c}):C\rightarrow X_{d}$ for a destination $d\in\mathcal{D}$ as 
\begin{align}
	\mathbf{F}_{d}(\mathbf{c})\equiv \arg. \min_{\mathbf{x}_{d}\in X_{d}}. Z_{d}(\mathbf{x}_{d}).
\end{align}
Then, the basic properties of the Legendre transformation lead to the following series of lemmas that depicts useful properties of the conjugate function~\citep[for a detailed explanation of the properties of Legendre transformation, see for example,][]{sandholm2010population}:

\begin{lemm}
	The optimal solution map $\mathbf{F}_{d}(\mathbf{c})$, $\forall d\in\mathcal{D}$ for a given link cost pattern $\mathbf{c}\in C$ is continuous and differentiable on $C$.
\end{lemm}

\begin{lemm}\label{Lemm:OptimalHandExpectedCost}
	The optimal value function $H_{d}^{*}(\mathbf{c})$, $\forall d\in\mathcal{D}$ of the problems~\eqref{Eq:CP-NGEV1} and \eqref{Eq:CP-NGEV2} is continuous, twice differentiable, and concave on $C$.
	It is given by
	\begin{align}
		H_{d}^{*}(\mathbf{c}) \equiv Z_{d}(\mathbf{F}_{d}(\mathbf{c}))= \sum_{i\in \mathcal{N}}q_{id}S_{id}(\mathbf{c}),\label{Eq:OptMap_MinECost}
	\end{align}
	where $S_{id}(\mathbf{c})$ is the function representing the expected minimum cost, which is defined in Eq.~\eqref{Eq:NGEV_ExpMinCost}.
\end{lemm}

\begin{lemm}\label{Lemm:ConjugateIdentity}
	The functions $H_{d}$, $H_{d}^{*}$, and $\mathbf{F}_{d}$, $\forall d\in\mathcal{D}$ satisfy the following identity:
	\begin{align}
		H_{d}^{*}(\mathbf{c}) + H_{d}(\mathbf{F}_{d}(\mathbf{c})) = \mathbf{c}\cdot \mathbf{F}_{d}(\mathbf{c}),
	\end{align}
	and the gradient of $H_{d}^{*}$ yields the optimal solution map $\mathbf{F}_{d}$:
	\begin{align}
		\mathbf{F}_{d}(\mathbf{c}) = \nabla H_{d}^{*}(\mathbf{c}).
	\end{align}
\end{lemm}

\begin{lemm}\label{Lemm:JacobMonotone}
	For $d\in\mathcal{D}$, the Jacobian matrix $\nabla \mathbf{F}_{d}(\mathbf{c}) = \nabla^{2}H_{d}^{*}(\mathbf{c})$ is symmetric and negative-semidefinite, and the negative of the optimal solution map, $-\mathbf{F}_{d}(\mathbf{c})$, is monotone on $C$:
	\begin{align}
		-(\mathbf{F}_{d}(\mathbf{c}') - \mathbf{F}_{d}(\mathbf{c}''))\cdot (\mathbf{c}' - \mathbf{c}'')\geq 0,\quad 
		\forall \mathbf{c}',\mathbf{c}''\in C.
	\end{align}
\end{lemm}

\subsection{Fixed-point problems}
The optimal solution map $\mathbf{F}_{d}$ defined above enables us to concisely represent the equilibrium conditions as a fixed-point problem:

\begin{prop}
	The NGEV equilibrium assignment is equivalent to the following fixed-point problem:
	\begin{align}
		\text{[FP-D]}
		\quad 
		\text{Find $\mathbf{x}^{*}\in X\cap Y$ such that}\quad 
		\mathbf{x}^{*}_{d} = \mathbf{F}_{d}(\mathbf{C}(\mathbf{x}^{*})),\quad \forall d\in\mathcal{D}.
	\end{align}
\end{prop}
\noindent This fixed-point problem above implies that, even if we know only the aggregated equilibrium link flows $\mathbf{x}^{*}$ or link costs $\mathbf{C}(\mathbf{x}^{*})$ (i.e., the arguments in the right-hand side), the disaggregated (or destination-specific) equilibrium link flows $\mathbf{x}_{d}$ can be obtained by evaluating the function $\mathbf{F}_{d}$.

This observation leads us to formulate the equilibrium conditions using only the aggregated link flows or link costs.
We define the NGEV-based link demand (total link flow) function $\mathbf{F}:C\rightarrow X$ as 
\begin{align}
	\mathbf{F}(\mathbf{c}) \equiv \sum_{d\in\mathcal{D}}\mathbf{F}_{d}(\mathbf{c}).
\end{align}
Then, we obtain the following fixed-point problems:
\begin{coro}
	The NGEV equilibrium assignment is equivalent to the following problems:
	\begin{align}
		&
		\text{[FP-A]}\quad \text{Find $\mathbf{x}^{*}\in X\cap Y$ such that}\quad 
		\mathbf{x}^{*} = \mathbf{F}(\mathbf{C}(\mathbf{x}^{*})),\\
		&
		\text{[FP-C]}\quad \text{Find $\mathbf{c}^{*}\in C$ such that}\quad 
		\mathbf{c}^{*}=\mathbf{C}(\mathbf{F}(\mathbf{c}^{*})).
	\end{align}
\end{coro}
We further obtain the equivalent system of non-linear equations as follows:
\begin{coro}\label{Prop:SNLE-C}
The NGEV equilibrium assignment is equivalent to the following system of non-linear equations:
\begin{align}
	\text{[NLE-C]}\quad 
	\text{Find $\mathbf{c}^{*}\in C$ such that}\quad \mathbf{C}^{-1}(\mathbf{c}^{*}) = \mathbf{F}(\mathbf{c}^{*}).\label{Eq:FP_Result}
\end{align}	
\end{coro}

\begin{remark}
The solution $\mathbf{x}^{*}$ of [FP-A] has a one-to-one correspondence with the solution $\mathbf{c}^{*}$ of [NLE-C] because the latter problem can be obtained from the former by a simple change of variables $\mathbf{c}^{*} = \mathbf{C}(\mathbf{x}^{*})$ and $\mathbf{x}^{*} = \mathbf{C}^{-1}$ based on the link cost function and its inverse defined in \textbf{Section~\ref{Sec:Settings}}.
It is also obvious that the solution $\mathbf{c}^{*}$ of [NLE-C] is exactly the same as that of [FP-C] because transforming the both sides of the former equation by the monotone link cost function $\mathbf{C}$ results in the latter equation, which implies the equivalency of the two equation. 

Also, the equation in [NLE-C] represents the demand-supply equilibrium condition: the link costs c are adjusted so as to equalize the ``supply-side link flows" (the left-hand side) and the ``demand-side link flows" (the right-hand side). 
For a general demand function, which is not necessarily based on the RU theory, it could be defined as the following non-linear complementarity problem:
\begin{align}
	\quad 
	\text{Find $\mathbf{c}^{*}\in C$ such that}\quad 
	\mathbf{0}\leq [\mathbf{c}^{*}-\mathbf{C}(\mathbf{0})] 
	\perp
	[\mathbf{C}^{-1}(\mathbf{c}^{*}) - \mathbf{F}(\mathbf{c}^{*})]\geq \mathbf{0}.\label{Eq:NLCP-SUE}
\end{align}
That is, in addition to the demand-supply equality condition (i.e. $C_{ij}^{-1}(\mathbf{c}^{*}) = F_{ij}(\mathbf{c}^{*})$ if $c_{ij}^{*} > C_{ij}(\mathbf{0})$), [NLE-C] also represents the possibility of an ``excess supply" case (i.e. $C_{ij}^{-1}(\mathbf{c}^{*})\geq F_{ij}(\mathbf{c}^{*}))$ if $c_{ij}^{*}=C_{ij}(\mathbf{0})$).
The inequality case, however, never occur in our model because the NGEV-based demand function $\mathbf{F}(\mathbf{c}^{*})$ is based on the RU theory, which implies that every alternative (link) is always chosen with a strictly positive (non-zero) probability (i.e. $\mathbf{c}^{*} = \mathbf{C}(\mathbf{x}) > \mathbf{C}(\mathbf{0})$ for some $\mathbf{x}>\mathbf{0}$).
Thus, [NLE-C] is sufficient for the description of the NGEV equilibrium assignment.
\end{remark}




\subsection{Variational inequality problems}
In addition to the fixed-point problem, the equilibrium conditions can be formulated as the variational inequality problems.
Specifically, the pair of functions $H_{d}$ and $H_{d}^{*}$ defined in the previous subsection leads to the following representation of the equilibrium conditions:
\begin{prop}\label{Prop:VI-SUE}
	The NGEV equilibrium assignment is equivalent to the following variational inequality problems:
	\begin{align}
		&\text{[VI-D]}
		&&\text{Find $\hat{\mathbf{x}}^{*}\in \hat{X}$ satisfying $\mathbf{x}^{*}\in X\cap Y$ such that }\notag\\
		&
		&&\quad \sum_{d\in\mathcal{D}}[\mathbf{C}(\mathbf{x}^{*}) - \nabla H_{d}(\mathbf{x}^{*}_{d})]\cdot (\mathbf{x}_{d} - \mathbf{x}^{*}_{d})\geq 0,
		\quad \forall \mathbf{x}\in X\cap Y,\\
		&\text{[VI-C]}
		&&\text{Find $\mathbf{c}^{*}\in C$ such that }\notag\\
		&
		&&\quad \sum_{d\in\mathcal{D}}[\mathbf{C}^{-1}(\mathbf{c}^{*}) - \nabla H^{*}_{d}(\mathbf{c}^{*})]\cdot (\mathbf{c} - \mathbf{c}^{*})\geq 0,
		\quad \forall \mathbf{c}\in C.
	\end{align}
\end{prop}

We present some remarks regarding these variational inequalities. 
First, the primal variational inequality [VI-D] implies that the vector field $-[\mathbf{C}(\mathbf{x}) - \nabla H_{d}(\mathbf{x}_{d})]$ is normal to $X_{d}$ for all $d\in\mathcal{D}$, at the equilibrium $\mathbf{x}^{*}$.
A similar characterization of the equilibrium has been exploited by \cite{Smith2016a}.
They extended the original Smith dynamic (route swapping dynamic) for deterministic equilibrium assignment in \cite{Smith1984} to the logit-based SUE assignment in a path flow space.

Second, the variational inequalities [VI-D] and [VI-C] are a natural primal-dual pair in the sense that a pair of ``inverse functions" for the demand-side ($\nabla H_{d}$ and $\nabla H_{d}^{*})$ and supply-side ($\mathbf{C}(\mathbf{x})$ and $\mathbf{C}^{-1}(\mathbf{c})$) is used.
This can be demonstrated by considering the integrable case, in which the variational inequalities in \textbf{Proposition~\ref{Prop:VI-SUE}} reduce to optimization problems, whereas our link cost function $\mathbf{C}$ or its inverse $\mathbf{C}^{-1}$ defined in \textbf{Section~\ref{Sec:Settings}} is not integrable.
If the link cost function is integrable, [VI-D] in \textbf{Proposition~\ref{Prop:VI-SUE}} is reduced to the following convex programming:
\begin{align}
	&\text{[CP-D]}
	&&\min_{\hat{\mathbf{x}}\in \hat{X}}. Z^{*}_{\text{SUE}}(\hat{\mathbf{x}}) 
	\equiv 
	\sum_{d\in\mathcal{D}}\oint_{\mathbf{x}}[\mathbf{C}(\mathbf{w})
	 - \nabla H_{d}(\mathbf{w})]\mathrm{d}\mathbf{w}\notag\\
	 &
	 &&
	 \hspace{20mm}= Z_{\text{UE}}(\mathbf{x}) - \sum_{d\in \mathcal{D}}H_{d}(\mathbf{x}_{d}),
\end{align}
where the function $Z_{\text{UE}}:\mathbb{R}_{+}^{|\mathcal{L}|}\rightarrow \mathbb{R}_{+}$ is defined as 
\begin{align}
	Z_{\text{UE}}(\mathbf{x})\equiv \oint_{\mathbf{x}}\mathbf{C}(\mathbf{w})\mathrm{d}\mathbf{w}.
\end{align}
This is a generalized version of the logit-based equivalent convex programming in \cite{Akamatsu1996-us,Akamatsu1997-bo}.
If we further assume that the inverse link cost function $\mathbf{C}^{-1}$ exists and is integrable, [VI-C] in \textbf{Proposition~\ref{Prop:VI-SUE}} reduces to the equivalent optimization problem in \cite{Daganzo1982-pd}:
\begin{align}
	&\text{[CP-C]}
	&&\min_{\mathbf{c}\geq \mathbf{0}}. Z^{*}_{\text{SUE}}(\mathbf{c}) \equiv 
	\sum_{d\in\mathcal{D}}\oint_{\mathbf{c}}[\mathbf{C}^{-1}(\mathbf{v}) - \nabla H_{d}^{*}(\mathbf{v})]\mathrm{d}\mathbf{v} \\
	&
	&&
	\hspace{20mm}= - Z^{*}_{\text{UE}}(\mathbf{c}) + \sum_{d\in \mathcal{D}}H_{d}^{*}(\mathbf{c}),
\end{align}
where the function $Z^{*}_{\text{UE}}$ is the Legendre transform of $Z_{\text{UE}}$:
\begin{align}
	Z^{*}_{\text{UE}}(\mathbf{c})\equiv 
	\max_{\mathbf{x}\geq \mathbf{0}}.[\mathbf{c}\cdot \mathbf{x} - Z_{\text{UE}}(\mathbf{x})]\equiv \oint_{\mathbf{c}}\mathbf{C}^{-1}(\mathbf{v})\mathrm{d}\mathbf{v}.
\end{align}
As shown in the next section, the pair of problems, [CP-D] and [CP-C], are useful for considering the stability of the day-to-day flow dynamic (even for the models with non-integrable link cost functions).

\begin{figure}[t]
	\centering
	\includegraphics[width=0.6\linewidth]{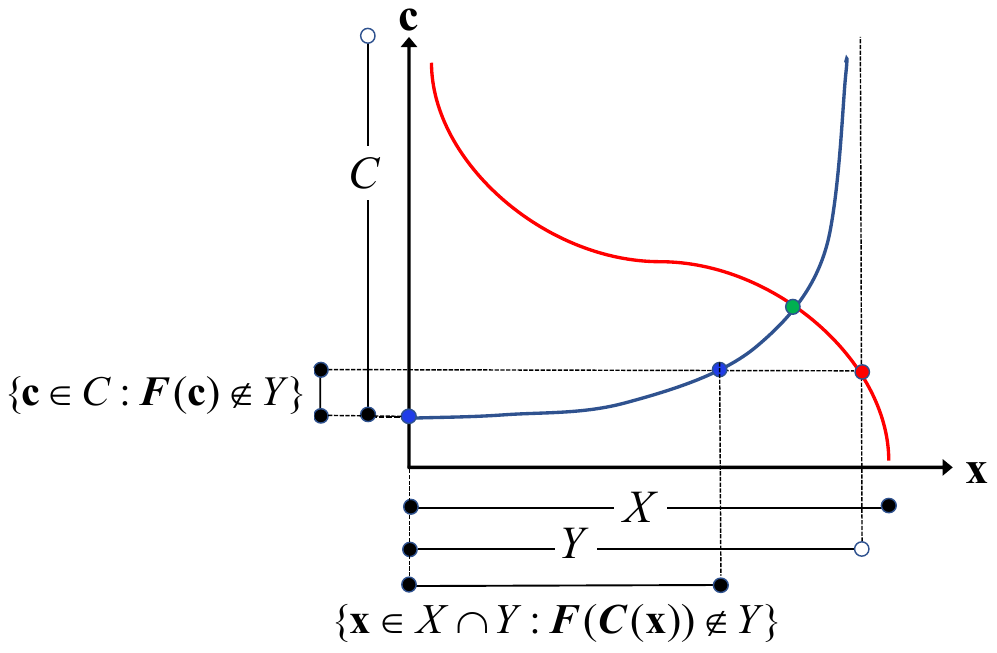}
	\caption{Illustrative example for $\mathbf{F}(\mathbf{C}(\mathbf{x}))\not\in Y$ in a simple 1-dimensional case}
	\label{Fig:Infeasibility}
\end{figure}

\section{Global stability of equilibrium}\label{Sec:ThoreticalAnalysis}
This section establishes the existence, uniqueness, and stability of the NGEV equilibrium.
These theoretical properties are usually proven in order from the existence by employing a standard theorem regarding variational inequality and fixed-point problems, such as Brouwer's fixed-point theorem.
However, such a standard theorem cannot easily be applied to our problem owing to the capacity constraints.
First, there may be the case that $\mathbf{F}(\mathbf{C}(\mathbf{x}))$ maps a feasible flow pattern $\mathbf{x}\in X\cap Y$ to a supply-infeasible flow pattern, i.e., $\mathbf{F}(\mathbf{C}(\mathbf{x}))\in X$ but $\mathbf{F}(\mathbf{C}(\mathbf{x}))\not\in Y$ (Figure~\ref{Fig:Infeasibility}).
This is because the condition $X\subseteq Y$ does not necessarily hold for our model (e.g., when OD demands are large relative to capacity limitations).
This complicates the application of the fixed-point theorem because the function $\mathbf{F}(\mathbf{C}(\mathbf{x}))$ does not map from a feasible region to itself.
Second, the feasible set is not a closed set owing to the strict inequality of the capacity constraints.


To address these problems, we prove the theoretical properties by an approach utilizing day-to-day dynamics and Lyapunov functions, referring to the works of \cite{smith1979existence} and \cite{Smith2015-my}.
We define day-to-day dynamics of link flows and costs, and present the Lyapunov functions for the dynamics.
Next, we prove the existence by investigating the evolution process of a traffic state under the day-to-day dynamics and showing the convergence to an equilibrium point within the feasible region.
Further, we prove the uniqueness by an approach that does not require the feasible set to be a closed set.
Then, we establish the global convergence of the day-to-day dynamics to the NGEV equilibrium using the Lyapunov stability theory.

\subsection{Day-to-day dynamics}
Here, we define the day-to-day dynamics of link flows and costs, which become the basis for the analysis of the theoretical properties.
Let $t\in\mathbb{R}_{+}$ be a day-to-day (continuous) time and $\mathbf{x}_{d}(t)$ be the disaggregated link flows at the time $t$.
We then consider the following day-to-day dynamic of link flows:
\begin{defi}[\textit{NGEV Dynamic-D}]\label{Defi:NGEVDynamicsD}
The NGEV day-to-day dynamic of the disaggregated link flows with the initial condition $\mathbf{x}(0)\in X\cap Y$ satisfying $\hat{\mathbf{x}}(0)\in\hat{X}$ is defined as
\begin{align}
\dot{\mathbf{x}}_{d}(t)=\mathbf{F}_{d}(\mathbf{C}(\mathbf{x}(t))) - \mathbf{x}_{d}(t),\quad \forall d\in\mathcal{D},\quad
\end{align}
where the dot denotes the derivative with respect to $t$ (i.e., $\mathrm{d}/\mathrm{d}t$).
\end{defi}



\noindent Under the day-to-day dynamic, excess link flows compared to the optimal solution map swap toward other links.
This dynamic is straightforwardly derived from the equivalent fixed-point problem [FP-A] formulated in the previous section.
Note that this dynamic can be interpreted as a link-based formulation of the logit-type evolutionary dynamic of a population game, which describes the evolution of strategy (i.e., route) distributions, in \cite{sandholm2010population}.

We also define the day-to-day dynamic of link costs as follows (hereinafter, we omit the argument $(t)$ for simplicity):
\begin{defi}[\textit{NGEV Dynamic-C}]\label{Defi:NGEVDynamicsC}
The NGEV day-to-day dynamic of link costs with the initial condition $\mathbf{c}(0)\in \mathbb{R}_{+}^{L}$ is defined by
\begin{align}
\dot{\mathbf{c}} = \mathbf{Q}[\mathbf{F}(\mathbf{c}) - \mathbf{C}^{-1}(\mathbf{c})],
\end{align}
where $\mathbf{Q}$ is a symmetric matrix, satisfies $\mathbf{Q}\mathbf{1}=\mathbf{0}$, and is positive definite\footnote{We introduce this matrix to satisfy the unit consistency of the equation: the unit of $\mathbf{F}(\mathbf{c}) - \mathbf{C}^{-1}(\mathbf{c})$ in the RHS is related to the \textit{flow}, and is thus converted into the unit related to the \textit{cost} by the matrix.}.
\end{defi}
\noindent We construct this dynamic based on the system of non-linear equations [NLE-C] in \textbf{Corollary~\ref{Prop:SNLE-C}} but not on the fixed-point problem [FP-C], unlike [NGEV Dynamic-D].
The reason why we do not use [FP-C] is that the mapping $\mathbf{C}(\mathbf{F}(\mathbf{c}))$ may not be appropriately defined for a given non-equilibrium cost pattern $\mathbf{c}$: $\mathbf{C}(\mathbf{F}(\mathbf{c}))$ cannot be calculated if $\mathbf{F}(\mathbf{c})\not\in Y$, as indicated at the beginning of this section.


\subsection{Lyapunov function for the dynamical system}
As we have seen in [VI-D] of \textbf{Proposition~\ref{Prop:VI-SUE}}, the equilibrium $\mathbf{x}^{*}$ can be characterized by $-[\mathbf{C}(\mathbf{x}^{*}) - \nabla H_{d}(\mathbf{x}^{*}_{d})]$, which is a vector field on $X\cap Y$, or indicates the ``downhill" direction of a field of ``force".
Therefore, the first step toward examining the stability property of the NGEV dynamic is determining whether the dynamic conforms with the vector field.
The following lemma provides a positive answer to this question: 

\begin{lemm}\label{Lemm:PositiveCorrelation}
The NGEV-dynamic of disaggregated link flows satisfies
\begin{align}
-\sum_{d\in\mathcal{D}}\left[ \mathbf{C}(\mathbf{x}) - \nabla H_{d}(\mathbf{x}_{d}) \right]\cdot \dot{\mathbf{x}}_{d}\geq 0,
\quad \forall \mathbf{x}\in X\cap Y,\forall \hat{\mathbf{x}}\in\hat{X}.
\end{align}
\end{lemm}
\begin{prf}
Theorem 6.2.10 (virtual payoff correlation) in \cite{sandholm2010population} can be applied.\qed
\end{prf}

Based on this lemma, we establish the following existence theorem of a \textit{Lyapunov function} for the day-to-day dynamic:
\begin{theo}[Existence of a Lyapunov function]\label{Theo:LyapuFunc}
Let $V$: $\hat{X} \rightarrow \mathbb{R}_{+}$ be a continuously differentiable function given by
\begin{align}
V(\hat{\mathbf{x}}) \equiv \mathbf{C}(\mathbf{x})\cdot \mathbf{x} - \sum_{d\in\mathcal{D}}\left[ H_{d}(\mathbf{x}_{d}) + H^{*}_{d}(\mathbf{C}(\mathbf{x})) \right].\label{Eq:LyapunovFunc}
\end{align}
Then, the function $V$ is a strict Lyapunov function for [NGEV Dynamic-D]: letting $\hat{\mathbf{x}}^{*}$ be an equilibrium state defined in [VI-D],

\begin{enumerate}
\renewcommand{\labelenumi}{(\roman{enumi})}
\item $V(\hat{\mathbf{x}}) > 0$ for all $\hat{\mathbf{x}}\in \hat{X}\setminus\{ \hat{\mathbf{x}}^{*} \}$,
\item $V(\hat{\mathbf{x}}) = 0$ if and only if $\hat{\mathbf{x}} = \hat{\mathbf{x}}^{*}$ and 
\item $\nabla V(\hat{\mathbf{x}}) < 0$ for all $\hat{\mathbf{x}}\in X\setminus\{ \hat{\mathbf{x}}^{*} \}$ and $\nabla V(\hat{\mathbf{x}}) = 0$ if and only if $\hat{\mathbf{x}} = \hat{\mathbf{x}}^{*}$.
\end{enumerate}
\end{theo}
\begin{prf}

(i) Consider a feasible flow pattern $\hat{\mathbf{x}}\neq \hat{\mathbf{x}}^{*}$.
It follows from the definition of $H_{d}^{*}(\mathbf{c})$ and $\mathbf{F}_{d}(\mathbf{c}(\mathbf{x}))$ (\textbf{Lemma~\ref{Lemm:OptimalHandExpectedCost}}) that the following inequality holds for an arbitrary feasible flow pattern $\mathbf{y}_{d}\neq \mathbf{F}_{d}(\mathbf{c}(\mathbf{x}))$:
\begin{align}
	H_{d}^{*}(\mathbf{C}(\mathbf{x}))\equiv \min_{\mathbf{x}'_{d}\in X_{d}}.[\mathbf{C}(\mathbf{x})\cdot \mathbf{x}'_{d} - H_{d}(\mathbf{x}'_{d})]<\mathbf{C}(\mathbf{x})\cdot \mathbf{y}_{d} - H_{d}(\mathbf{y}_{d}),\quad \forall d\in\mathcal{D}.\label{Eq:HOptimalandH}
\end{align}
Hence, for any feasible flow patterns $\hat{\mathbf{x}}\neq \hat{\mathbf{x}}^{*}$ and $\mathbf{y}\neq \mathbf{F}(\mathbf{C}(\mathbf{x}))$, the function $V$ always satisfies
\begin{align}
	V(\hat{\mathbf{x}})
	&\equiv \mathbf{C}(\mathbf{x})\cdot \mathbf{x}
	- \sum_{d\in\mathcal{D}}[H_{d}(\mathbf{x}_{d}) + H_{d}^{*}(\mathbf{C}(\mathbf{x}))]\\
	&>\mathbf{C}(\mathbf{x})\cdot \mathbf{x}
	- \sum_{d\in\mathcal{D}}[H_{d}(\mathbf{x}_{d}) + \mathbf{C}(\mathbf{x})\cdot \mathbf{y}_{d} -H_{d}(\mathbf{y}_{d})]\\
	&=\mathbf{C}(\mathbf{x})\cdot (\mathbf{x}-\mathbf{y})
	- \sum_{d\in\mathcal{D}}[H_{d}(\mathbf{x}_{d}) -H_{d}(\mathbf{y}_{d})].
\end{align}
Setting $\hat{\mathbf{x}} = \hat{\mathbf{y}}\equiv\{\mathbf{y}_{1},\ldots, \mathbf{y}_{D} \}$ in the right-hand side of the inequality above, we obtain $V(\hat{\mathbf{x}})>0$ for all feasible flow patterns that are not equilibrium.\qed

(ii) As $H_{d}^{*}$ is a Legendre transformation of $H_{d}$, the following equality always holds (\textbf{Lemma~\ref{Lemm:ConjugateIdentity}}):
\begin{align}
	&
	H_{d}^{*}(\mathbf{C}(\mathbf{x})) + H_{d}(\mathbf{F}_{d}(\mathbf{C}(\mathbf{x}))) = \mathbf{C}(\mathbf{x})\cdot \mathbf{F}_{d}(\mathbf{C}(\mathbf{x})),
	&&
	\forall x\in X\cap Y.
\end{align}
Substituting the definition $\mathbf{x}_{d}^{*} = \mathbf{F}_{d}(\mathbf{C}(\mathbf{x}))$ of the equilibrium $\mathbf{x}_{d}^{*}$ into this equality, we obtain
\begin{align}
	&
	H_{d}^{*}(\mathbf{C}(\mathbf{x}^{*})) + H_{d}(\mathbf{x}^{*}) = \mathbf{C}(\mathbf{x}^{*})\cdot \mathbf{x}_{d}^{*}
	&&
	\forall d\in\mathcal{D},
\end{align}
which implies that $V(\hat{\mathbf{x}}^{*}) = \mathbf{C}(\mathbf{x}^{*})\cdot \mathbf{x}^{*} - \mathbf{C}(\mathbf{x}^{*})\cdot \mathbf{x}^{*}=0$ if $\hat{\mathbf{x}}^{*}$ is equilibrium.
The ``only if" part is obvious from (i).\qed

(iii) A straightforward calculation of $\dot{V}(\hat{\mathbf{x}})$ yields
\begin{align}
	\mathrm{d}V(\hat{\mathbf{x}})/\mathrm{d}t 
	&= \nabla_{\mathbf{x}}V\cdot \mathrm{d}\mathbf{x}/\mathrm{d}t + \nabla_{\hat{\mathbf{x}}}V\cdot \mathrm{d}\hat{\mathbf{x}}/\mathrm{d}t\\
	&=[\mathbf{C}(\mathbf{x}) + \nabla_{\mathbf{x}}\mathbf{C}(\mathbf{x})\cdot \mathbf{x}
	- \nabla_{\mathbf{x}}\mathbf{C}(\mathbf{x})\cdot \mathbf{F}(\mathbf{C}(\mathbf{x}))]\cdot \mathbf{x}
	- \sum_{d\in\mathcal{D}}\nabla_{\mathbf{x}_{d}}H_{d}(\mathbf{x}_{d})\cdot \dot{\mathbf{x}}_{d}\\
	&=[\mathbf{C}(\mathbf{x}) - \nabla_{\mathbf{x}}\mathbf{C}(\mathbf{x})\cdot (\mathbf{F}(\mathbf{C}(\mathbf{x})) - \mathbf{x})]\cdot \dot{\mathbf{x}}
	- \sum_{d\in\mathcal{D}}\nabla_{\mathbf{x}_{d}}H_{d}(\mathbf{x}_{d})\cdot \dot{\mathbf{x}}_{d}\\
	&=-\dot{\mathbf{x}}^{\mathrm{T}}[\nabla_{\mathbf{x}}\mathbf{C}(\mathbf{x})]\dot{\mathbf{x}} + \sum_{d\in\mathcal{D}}[\mathbf{C}(\mathbf{x})-\nabla_{\mathbf{x}_{d}}H_{d}(\mathbf{x}_{d})]\cdot \dot{\mathbf{x}}_{d}.
\end{align} 
The first term of $\dot{V}$ is negative for all $\dot{\mathbf{x}}\neq \mathbf{0}$ because $\mathbf{C}(\mathbf{x})$ is strictly monotone (\textbf{Proposition~\ref{Prop:CostMonotone}}), which implies that the Jacobian $\nabla \mathbf{C}(\mathbf{x})$ is positive definite.
As the second term is also negative for all $\dot{\mathbf{x}}\neq \mathbf{0}$ (\textbf{Lemma~\ref{Lemm:PositiveCorrelation}}), $\dot{V}\leq 0$ always holds, where the equality holds if and only if $\dot{\mathbf{x}} = 0$ at the equilibrium $\mathbf{x}^{*}$.\qed
\end{prf}

For the dynamic [NGEV Dynamic-C] of link costs, we obtain the following Lyapunov function in a straightforward manner:
\begin{theo}[Lyapunov function for the link cost dynamic]\label{Theo:LyapuFunc-C}
Let the inverse cost function $\mathbf{C}^{-1}:C\rightarrow Y$ be given by \textbf{Proposition~\ref{Prop:InverseCost}}, and let the link cost dynamic be given by [NGEV Dynamic-C].
Define the function $V_{C}:C\rightarrow \mathbb{R}_{+}$ as
\begin{align}
	V_{C}(\mathbf{c})\equiv 
    \cfrac{1}{2}\| \mathbf{Q} [\mathbf{F}(\mathbf{c})-\mathbf{C}^{-1}(\mathbf{c})]\|^{2}.\label{Eq:Lyapunov_Cost}
\end{align}
	Then, the function is a strict Lyapunov function for the dynamic, that is, $V_{C}(\mathbf{c}(t))$ is strictly decreasing, away from the equilibrium, for any solution $\mathbf{c}(\cdot)$ and attains its minimum at the solution $\mathbf{c}^{*}$ of [FP-C]:
	\begin{align}
		\min_{\mathbf{c}\in C}V_{C}(\mathbf{c})=V_{C}(\mathbf{c}^{*}) = 0.
	\end{align}
\end{theo}
\begin{prf}
Because the function $V_{C}$ is quadratic, it is obvious that (i) $V_{C}(\mathbf{c})\geq 0$ for all $\mathbf{c}\in C$, (ii) $V_{C}(\mathbf{c}) = 0$ if and only if $\mathbf{c}$ is equilibrium.
We thus prove that $\dot{V}_{C}(\mathbf{c})<0$ if $\mathbf{c}$ is not equilibrium.
A simple calculation derives the following:
\begin{align}
	\mathrm{d}V_{C}(\mathbf{c})/\mathrm{d}t 
	&= \nabla_{\mathbf{c}}V_{C}\cdot \mathrm{d}\mathbf{c}/\mathrm{d}t\\
	&=\mathbf{Q}[\nabla \mathbf{F}(\mathbf{c}) - \nabla \mathbf{C}^{-1}(\mathbf{c})]\cdot (\mathbf{Q}[\mathbf{F}(\mathbf{c}) - \mathbf{C}^{-1}(\mathbf{c})])\cdot \dot{\mathbf{c}}\\
	&=\dot{\mathbf{c}}^{T}\mathbf{M}\dot{\mathbf{c}},
\end{align}
where the matrix $\mathbf{M}$ is defined as 
\begin{align}
	 \mathbf{M}\equiv \mathbf{Q}[\nabla \mathbf{F}(\mathbf{c}) - \nabla \mathbf{C}^{-1}(\mathbf{c})].
\end{align}
\noindent To prove that $\dot{V}_{C}(\mathbf{c})<0$, it is sufficient to prove that the matrix $\mathbf{M}$ is strictly negative-definite.

First, the matrix $\nabla \mathbf{F}(\mathbf{c})$ is semi-negative definite because $-\mathbf{F}(\mathbf{c}) = -\sum_{d\in\mathcal{D}}\mathbf{F}_{d}(\mathbf{c})$ is monotone (\textbf{Lemma~\ref{Lemm:JacobMonotone}}).
The matrix $\nabla \mathbf{C}^{-1}(\mathbf{c})$ is positive definite because $\mathbf{C}^{-1}(\mathbf{c})$ is strictly monotone (\textbf{Proposition~\ref{Prop:InverseCost}}).
It follows that the matrix $[\nabla \mathbf{F}(\mathbf{c})-\nabla \mathbf{C}^{-1}(\mathbf{c})]$ is strictly negative definite, and the matrix $\mathbf{M}$ is the product of the positive-definite and negative definite matrices.
By combining this fact and Hine's lemma~\citep[see \textbf{Lemma 8.5.1} in][]{sandholm2010population}, we show that each eigenvalue of the matrix $\mathbf{M}$ has a negative real part.
Hence, the matrix $\mathbf{M}$ is strictly negative definite, which implies that $\dot{V}_{C}(\mathbf{c})\leq 0$ for all $\mathbf{c}\in C$, where the equality holds if and only if $\dot{\mathbf{c}} = 0$ at the equilibrium $\mathbf{c}^{*}$.\qed

\end{prf}

\subsection{Some remarks on the Lyapunov function}
First, our Lyapunov theorems are extended applications of Theorem 5.2 in \cite{hofbauer2009stable} (or Theorem 7.2.8 in \cite{sandholm2010population}) in terms of having a non-simplex and capacity-constrained solution space; this implies that the theorems are applicable not only to the NGEV equilibrium but also the entire perturbed equilibrium reproduced by a route choice model consistent with the RU theory.
\cite{hofbauer2009stable} showed that the alternative or its choice probability distribution consistent with an arbitrary RU model is derived by the Legendre transform of \textit{admissible deterministic perturbation}, which is a function representing the stochastic effect.
Such admissible deterministic perturbation (and the conjugate function) corresponds to $H_{d}$ (and $H^{*}_{d}$) in this study. and it is formulated as an entropy function in response to the postulated NGEV route choice model.
Meanwhile, the Lyapunov theorem shows that once these functions are obtained, we can construct day-to-day dynamics with the virtual positive correlation and establish the Lyapunov function for the dynamics.
This suggests that we can establish the global stability in a constructive manner for an arbitrary RU model.
To the best of our knowledge, this is the first application of this theorem to traffic assignment problems.

Next, to provide an intuitive meaning of the Lyapunov function $V$, we consider cases where link cost functions are integrable.
From the definition of $Z_{\text{UE}}$ and $Z_{\text{UE}}^{*}$ in equivalent convex programming, [CP-D] and [CP-C], it follows that 
\begin{align}
	\mathbf{C}(\mathbf{x})\cdot \mathbf{x} = 
	Z_{\text{UE}}(\mathbf{x}) + Z_{\text{UE}}^{*}(\mathbf{C}(\mathbf{x}))
\end{align}
always hold for all $\mathbf{x}\in Y$.
Substituting this into the definition of $V$, we have 
\begin{align}
	V(\hat{\mathbf{x}}) 
	&= [Z_{\text{UE}}(\mathbf{x}) + Z_{\text{UE}}^{*}(\mathbf{C}(\mathbf{x}))] 
	- \sum_{d\in\mathcal{D}}[H_{d}(\mathbf{x}_{d}) + H_{d}^{*}(\mathbf{C}(\mathbf{x}))]\\
	&= [Z_{\text{UE}} - \sum_{d\in\mathcal{D}}H_{d}(\mathbf{x}_{d}) ]
	- [-Z_{\text{UE}}^{*}(\mathbf{C}(\mathbf{x})) + \sum_{d\in\mathcal{D}}H_{d}^{*}]\\
	&= Z_{\text{SUE}}(\mathbf{x}) - Z^{*}_{\text{SUE}}(\mathbf{C}(\mathbf{x})).
\end{align}
That is, the Lyapunov function $V$ represents the difference between the objective functions, $Z_{\text{SUE}}$ and $Z_{\text{SUE}}^{*}$, of the primal and dual convex programming.
The duality theorem for convex programming shows that $Z_{\text{SUE}}(\mathbf{x})\geq Z_{\text{SUE}}^{*}(\mathbf{C}(\mathbf{x}))$ (i.e., $V(\hat{\mathbf{x}}) \geq 0$) for all feasible flow patterns $\mathbf{x}$, and $Z_{\text{SUE}}(\mathbf{x}) = Z_{\text{SUE}}^{*}(\mathbf{C}(\mathbf{x}))$ (i.e., $V(\hat{\mathbf{x}}) = 0$) if and only if $\mathbf{x}$ is equilibrium.

\subsection{Existence, uniqueness, and stability}
We are now ready to establish the theoretical properties of the NGEV equilibrium.
First, we establish the existence of the equilibrium by utilizing the NGEV Dynamic-D in \textbf{Definition~\ref{Defi:NGEVDynamicsD}} as follows:
\begin{theo}
Assume that $X\cap Y\neq \emptyset$.
Then, there exists an equilibrium state $\hat{\mathbf{x}}^{*}\in \hat{X}$ satisfying $X \cap Y$.
\end{theo}
\begin{prf}
Sketch of the proof is shown here.
We consider the trajectory of a traffic state under the day-to-day dynamic [NGEV Dynamic-D].
First, the property (iii) in \textbf{Theorem~\ref{Theo:LyapuFunc}} shows that for an arbitrary time $t\in\mathbb{R}_{+}$, the traffic state $\hat{\mathbf{x}}(t)$ belongs to the demand- and supply-feasible region; that is, $\hat{\mathbf{x}}(t)\in \hat{X}$ and $\mathbf{x}(t)\in X\cap Y$.
By utilizing this fact, it follows that a limit point $\hat{\mathbf{x}}^{\infty}$ of the dynamical system must belong to the feasible region and satisfy the following relationship: $V(\hat{\mathbf{x}}^{\infty}) = 0$.
This means that $\hat{\mathbf{x}}^{\infty}$ is an equilibrium state defined in [VI-D], and thus the equilibrium state belongs to the feasible region.
Therefore, an equilibrium state exists in the feasible region.
See \ref{Sec:App-Existence} for further details.\qed
\end{prf}

We next establish the uniqueness of the NGEV equilibrium as follows:
\begin{theo}
Assume that $X\cap Y\neq \emptyset$.
Then, the equilibrium $\hat{\mathbf{x}}^{*}$ of the solution to [VI-D] is unique.
\end{theo}
\begin{prf}
We use proof by contradiction. We assume that two different equilibrium states exist and show that this assumption leads to a contradiction.
See \ref{Sec:App-Unique} for further details.\qed
\end{prf}

We finally establish the stability of the NGEV equilibrium under the NGEV Dynamic-D as follows:
\begin{theo}
Assume that $X\cap Y\neq \emptyset$.
Then, the NGEV equilibrium under NGEV Dynamic-D is globally asymptotically stable.
\end{theo}
\begin{prf}
The application of the Lyapunov's stability theory to the function $V$ in \textbf{Theorem~\ref{Theo:LyapuFunc}} proves this theorem.\qed
\end{prf}
\noindent The stability under the NGEV Dynamic-C is also established from the Lyapunv function~\eqref{Eq:Lyapunov_Cost}. 

The stability result suggests that the proposed dynamics and Lyapunov functions naturally lead to the development of efficient solution algorithms, whose convergence is ensured.
Although these algorithms are not directly related to the contribution of this study (i.e., modeling and theoretical analysis), we introduce them in \ref{Sec:Algorithm} for interested readers.

\section{Concluding remarks}
This study presented an analytical framework for schedule-based transit assignment that describes the dynamics of public transport congestion and boarding queues at stations. 
To explicitly consider vehicle capacity constraints, we modeled the priority of passengers already on board to those boarding from a station by introducing boarding-queue links to the time-space network. The flow of a boarding-queue link is limited by the residual capacity, which is the link capacity subtracted by the flow of on-board passengers. 
Under this asymmetric link cost structure and explicit capacity constraints, the NGEV equilibrium assignment \citep{oyama2022markovian}, which captures the correlation among time-space paths without path enumeration, was formulated.
For theoretical analysis, the day-to-day dynamics of both link flows and link costs were derived. We showed that Lyapunov functions exist for the NGEV dynamics, with which the existence, uniqueness, and global stability of equilibrium were established.


As concluding remarks, we here want to emphasize that the novelty and contributions of this study are significant in both contexts of dynamic transit assignment and general traffic assignment.

First, in the literature, most dynamic transit assignment models with explicit capacity constraints were based on the deterministic equilibrium framework and relied on simulation.
Thus, the theoretical properties of the solution were not clarified. 
This study filled the gap by presenting an analytical framework for Markovian dynamic transit assignment. The guarantee of desirable solution properties is important for practical application in obtaining a reliable solution for dynamic transit assignment. The framework also describes doubly (i.e., within-day and day-to-day) dynamics of public transport congestion, which has often been neglected in the literature \citep{Nuzzolo2001transit, cats2016dynamic}. Moreover, the present theoretical analysis naturally leads to the development of efficient solution algorithms, whereas the method of successive averages, an iterative method with a poor convergence rate, has been widely used in dynamic transit assignment.

Second, we theoretically analyzed the stochastic equilibrium with explicit capacity constraints and an asymmetric link cost structure, which, to the best of our knowledge, has never been studied in the field of general traffic assignment. 
The present Lyapunov theorems are an extension of Theorem 5.2 of \cite{Hofbauer2002a} in terms of having a non-simplex and capacity-constrained solution space, and are applicable to equilibrium with any additive RU model. 
Moreover, because the proposed day-to-day dynamics are defined based on link-based variables, they do not suffer from the problems of initial path flow patterns and path overlapping \citep{He2010}. Thus, our theoretical analysis is an extension/generalization of existing studies on general traffic dynamics.

Possible future work includes applications to real networks and policy analysis, such as pricing or timetable design. For practical application, parameter estimation with real observations is required.
Moreover, the utilization of various emerging data sources would be an appealing direction.







\section*{Acknowledgements}
This work was financially supported by JSPS KAKENHI Grant numbers JP21H01448 and JP20H00265.


\appendix

\renewcommand{\thelemm}{\Alph{section}\arabic{lemm}}

\section{Proof of Proposition~\ref{Prop:VI-SUE}}
As a preliminary analysis, we establish the relationship between the NGEV equilibrium assignment and a complementarity problem.
First, the flow conservation condition can be expressed in the following complementarity form:
\begin{align}
&
\begin{cases}
	s_{id}\left[ \sum_{n\in\mathcal{F}(i)}x_{in}^{d} - \sum_{n\in\mathcal{B}(i)}x_{ni}^{d}-q_{i}^{d}\right]=0\\
	s_{id}\geq 0,
	\quad \sum_{n\in\mathcal{F}(i)}x_{in}^{d} - \sum_{n\in\mathcal{B}(i)}x_{ni}^{d}-q_{i}^{d}\geq 0,
\end{cases}
&&
\forall i\in\mathcal{N},\forall d\in\mathcal{D}.
\end{align}
This is revised in a matrix form as follows:
\begin{align}
	&
	\mathbf{0}\leq \left[\mathbf{A}_{d}\mathbf{x}_{d} - \mathbf{q}_{d}\right]\perp \mathbf{s}_{d} \geq \mathbf{0},
	&&
	\forall d\in\mathcal{D}.\label{Eq:App-PropVI_FCNode}
\end{align}

We next consider the link flow distribution condition.
In stochastic assignment based on a RU model, each alternative (i.e. link) is always chosen with a positive probability.
From this fact, we derive the following identity equivalent to Eq.~\eqref{Eq:LinkChoiceProb}: 
\begin{align}
	&
	\ln p_{ij\mid i}^{d} = 
	\ln \alpha_{ji}^{d}e^{-\theta_{i}^{d}(c_{ij}+s_{jd})}
	- \ln \sum_{n\in\mathcal{F}(i)}\alpha_{ni}^{d}e^{-\theta_{i}^{d}(c_{in}+s_{nd})}
	&&\forall ij\in\mathcal{L},\forall d\in\mathcal{D}.\label{Eq:App-PropVI_1}
\end{align}
From Eq.~\eqref{Eq:LinkChoiceProbSum}, we also have
\begin{align}
	&
	1 = e^{\theta_{i}^{d}s_{id}}\sum_{n\in\mathcal{F}(i)}\alpha_{ni}^{d}e^{-\theta_{i}^{d}(c_{in} + s_{nd})}
	&&\forall i\in\mathcal{N},\forall d\in\mathcal{D}.\label{Eq:App-PropVI_2}
\end{align}
Substituting Eq.~\eqref{Eq:App-PropVI_2} into Eq.~\eqref{Eq:App-PropVI_1} yoelds the following equation:
\begin{align}
	&\ln p_{ij\mid i}^{d} = -\theta_{i}^{d}(c_{ij} + s_{jd} - s_{id}) + \ln \alpha_{ji}^{d},
	&&\forall ij\in\mathcal{L},\forall d\in\mathcal{D}.
\end{align}
By utilizing these equations, the link flow distribution condition~\eqref{Eq:LFDist} can be expressed as follows:
\begin{align}
	&c_{ij} + \cfrac{1}{\theta_{i}^{d}}\ln \cfrac{x_{ij}^{d}}{\alpha_{ji}^{d}\sum_{n\in\mathcal{F}(i)}x_{in}^{d}} + s_{jd} - s_{id} = 0,
	&&
	\forall ij\in\mathcal{L},\forall d\in\mathcal{D}.\label{Eq:App-PropVI_3}
\end{align}
We then derive the equivalent complementarity condition as follows:
\begin{align}
&
\begin{cases}
	x_{ij}^{d}\left[	c_{ij} + \cfrac{1}{\theta_{i}^{d}}\ln \cfrac{x_{ij}^{d}}{\alpha_{ji}^{d}\sum_{n\in\mathcal{F}(i)}x_{in}^{d}} + s_{jd} - s_{id}	\right] = 0\\
	x_{ij}^{d}\geq 0,\quad c_{ij} + \cfrac{1}{\theta_{i}^{d}}\ln \cfrac{x_{ij}^{d}}{\alpha_{ji}^{d}\sum_{n\in\mathcal{F}(i)}x_{in}^{d}} + s_{jd} - s_{id} \geq 0,
\end{cases}
&&
\forall ij\in\mathcal{L},\forall d\in\mathcal{D}.
\end{align}
In addition, the partial derivative of the entropy function $H_{d}$ is given as follows:
\begin{align}
	&
	\cfrac{\partial H_{d}(\mathbf{x}^{d})}{\partial x_{ij}^{d}} = 
	-\cfrac{1}{\theta_{i}^{d}}
	\ln \cfrac{x_{ij}^{d}}{\alpha_{ji}^{d}\sum_{n\in\mathcal{F}(i)}x_{in}^{d}},
	&&
	\forall ij\in\mathcal{L},\forall d\in\mathcal{D}.
\end{align}
Thus, the complementarity condition is expressed in the following matrix form:
\begin{align}
	&
	\mathbf{0}\leq \left[\mathbf{C}(\mathbf{x}) - \nabla H_{d}(\mathbf{x}_{d}) - 
	\left(\mathbf{A}_{d}\right)^{\mathrm{T}}\mathbf{s}_{d} \right]
	\perp
	 \mathbf{x}_{d} \geq \mathbf{0},
	&&\forall d\in\mathcal{D}.\label{Eq:App-PropVI_LFDist}
\end{align}
It follows that the NGEV equilibrium assignment is equivalent to the complementarity problem, which determines $\langle \mathbf{x}^{*},\mathbf{s}^{*} \rangle\in \mathbb{R}_{+}^{|\mathcal{L}|\times (|\mathcal{N}|\times|\mathcal{D}|)}\cap Y$ satisfying Eqs.~\eqref{Eq:App-PropVI_FCNode} and \eqref{Eq:App-PropVI_LFDist}.
In addition, because we consider non-negative variables, the abovementioned complementarity problem is equivalent to the following variational inequality problem:
\begin{align}
	&\sum_{d\in\mathcal{D}}
	\left\{
	\left[\mathbf{A}_{d}\mathbf{x}^{*}_{d} - \mathbf{q}_{d}\right]
	\cdot 
	[\mathbf{s}_{d} - \mathbf{s}_{d}^{*}]
	+ \left[\mathbf{c}(\mathbf{x}^{*}) - \nabla H_{d}(\mathbf{x}^{*}_{d}) - 
	\left(\mathbf{A}_{d}\right)^{\mathrm{T}}\mathbf{s}_{d} \right]
	\cdot 
	[\mathbf{x}_{d} - \mathbf{x}^{*}_{d}]
	\right\}\geq 0,
	&&
	\forall (\mathbf{x},\mathbf{s})\in \mathbb{R}_{+}^{|\mathcal{L}|\times (|\mathcal{N}|\times|\mathcal{D}|)}\cap Y.\label{Eq:App-VINGEV}
\end{align}

We now establish the relationship between the NGEV equilibrium assignment and the variational inequality problems.
First, by setting the feasible region of $\mathbf{x}$ to the convex region $X$ satisfying Eq.~\eqref{Eq:FCNode}, the first term in Eq.~\eqref{Eq:App-VINGEV} becomes zero.
In addition, we have the following equation for the second term:
\begin{align}
	&\left(\mathbf{A}_{d}\right)^{\mathrm{T}}\mathbf{s}_{d}\cdot 
	[\mathbf{x}_{d} - \mathbf{x}^{*}_{d}]
	=(\mathbf{s}_{d})^{\mathrm{T}}
	\left[\mathbf{A}_{d}\mathbf{x}_{d} - \mathbf{A}_{d}\mathbf{x}^{*}_{d}\right]
	= 0,
	&&\forall d\in\mathcal{D}.
\end{align}
In summary, the NGEV equilibrium assignment is equivalent to the following variational inequality problem, which determines $\hat{\mathbf{x}}\in \hat{X}$ satisfying $\mathbf{x}^{*}\in X\cap Y$ such that
\begin{align}
	\sum_{d\in\mathcal{D}}
	\left[\mathbf{C}(\mathbf{x}^{*}) - \nabla H_{d}(\mathbf{x}^{*}_{d}) \right]
	\cdot 
	[\mathbf{x}_{d} - \mathbf{x}^{*}_{d}]\geq 0.
\end{align}
This completes the proof.

We here omit the proof of the dual of the variational inequality since the transformation of the NGEV equilibrium assignment to Eq.~\eqref{Eq:FP_Result} becomes the proof as follows.
Specifically, the dual variational inequality [VI-C] can be expressed as 
\begin{align}
	&
	[\mathbf{C}^{-1}(\mathbf{c}^{*}) - \mathbf{F}(\mathbf{c}^{*})]\cdot (\mathbf{c} - \mathbf{c}^{*})\geq 0,
	&&
	\forall \mathbf{c}\in C.\label{Eq:FP_Base}
\end{align}
It follows from this that [VI-C] is equivalent to Eq.~\eqref{Eq:FP_Result}, which is the equilibrium condition itself.


\qed

\section{Proof of the existence}\label{Sec:App-Existence}
For a given initial state $\hat{\mathbf{x}}^{0}\in \hat{X}$ satisfying $\sum_{d\in\mathcal{D}}\mathbf{x}_{d}^{0}\in X\cap Y$, we consider the trajectory of a traffic state $\hat{\mathbf{x}}(t)$ under the day-to-day dynamic [NGEV Dynamic-D].
We first show that the traffic state does not leave the feasible region as follows: 
\begin{lemm}\label{Lemm:ExistLemma1}
Consider a traffic state $\hat{\mathbf{x}}(t)$, $\forall t\in\mathbb{R}_{+}$ resulting from the day-to-day dynamic [NGEV Dynamic-D].
Then, the traffic state belongs to the demand- and supply-feasible region, that is,
\begin{align}
\hat{\mathbf{x}}\in\hat{X},\quad \text{and}
\quad 
\mathbf{x}(t)\in X\cap Y.
\end{align}
\end{lemm}
\begin{prf}
We first observe that the traffic state does not leave the region $\hat{X}$ (and $X$) from the following relationship:
\begin{align}
\mathbf{A}_{d}\cfrac{\mathrm{d}\mathbf{x}_{d}(t)}{\mathrm{d}t}= \mathbf{A}_{d}\mathbf{F}_{d}(\mathbf{C}(\mathbf{x})) - \mathbf{A}_{d}\mathbf{x}_{d} = \mathbf{0},\quad \forall d\in\mathcal{D}.
\end{align}
This means that for any $d\in\mathcal{D}$, the demand-feasible condition $\mathbf{A}_{d}\mathbf{x}_{d} = \mathbf{q}_{d}$ (given constant) must not be violated.
Hence, the disaggregated link flows belong to the region $\hat{X}$ and aggregated link flows belong to $X$.

We next show that the traffic state does not leave the region $Y$ by contradiction.
Suppose that the traffic state leaves the region $Y$ at $t = \tau\in\mathbb{R}_{+}$.
This implies that as $t$ approaches $\tau$, the state approaches the boundary of $Y$, which is denoted by $\mathrm{bd}(Y)$.
By definition, this also means there exist links where the link costs approach infinity, that is,
\begin{align}
\lim_{t \to \tau}C_{ij}(\mathbf{x}(t)) = \infty,\quad \exists ij\in\mathcal{L}.
\end{align}
We denote by $\mathcal{L}^{\infty}$ the set of such links whose costs approach infinity.

By considering Eq.~\eqref{Eq:HOptimalandH}, the function $V$ satisfies the following relationship for any feasible flow pattern $\mathbf{y}$:
\begin{align}
	V(\hat{\mathbf{x}}(t))
	 \geq 
	\sum_{ij\in\mathcal{L}^{\infty}}C_{ij}(\mathbf{x}(t))\cdot (x_{ij}(t)-y_{ij})
	+ \sum_{ij\in\mathcal{L}\setminus \mathcal{L}^{\infty}}C_{ij}(\mathbf{x}(t))\cdot (x_{ij}(t)-y_{ij})
	- \sum_{d\in\mathcal{D}}[H_{d}(\mathbf{x}_{d}(t)) -H_{d}(\mathbf{y}_{d})].
\end{align}
By setting $\hat{\mathbf{y}} = \hat{\mathbf{x}}(0)$, the first term obviously approaches infinity as $t$ approaches $\tau$, the second term approaches a finite value, and the third term approaches a finite value because the entropy function is bounded.
Thus, $\lim_{t \to \tau} V(\hat{\mathbf{x}}(t)) = \infty$.
However, this contradicts the property (iii) in \textbf{Theorem~\ref{Theo:LyapuFunc}}, where $V$ does not increase as time increases.
Hence, the traffic state does not leave $Y$.\qed
\end{prf}

Based on this lemma, we show the property of the limit point of this day-to-day dynamic as follows:
\begin{lemm}
Let $\hat{\mathbf{x}}^{\infty}$ be the limit point of the day-to-day dynamic [NGEV Dynamic-D].
The limit point belongs to the demand- and supply-feasible region, that is,
\begin{align}
\hat{\mathbf{x}}\in \hat{X},\quad \text{and}\quad 
\mathbf{x}^{\infty}\in X\cap Y.
\end{align}
Moreover, $V(\hat{\mathbf{x}}^{\infty}) = 0$.
\end{lemm}
\begin{prf}
Because the function $V$ is a monotonically decreasing function of $t$ and bounded below, there exists a limit point $\hat{\mathbf{x}}^{\infty}$ where $V$ is zero (the lower bound).
Next, as $\hat{\mathbf{x}}(t)\in\hat{X}$ and $\mathbf{x}(t)\in X\cap Y$, $\forall t\in\mathbb{R}_{+}$ (\textbf{Lemma~\ref{Lemm:ExistLemma1}}), the limit point belongs to the closure of the feasible region, that is,
\begin{align}
\hat{\mathbf{x}}\in \mathrm{cl}(\hat{X}) = \hat{X},\quad \text{and}\quad 
\mathbf{x}^{\infty}\in \mathrm{cl}(X\cap Y) = X\cap \mathrm{cl}(Y),
\end{align}
where $\mathrm{cl}$ denotes the closure of a set and satisfies the following relationship: $\mathrm{cl}(Y) = Y\cap \mathrm{bd}(Y)$.
We see that $\mathbf{x}^{\infty}$ cannot belong to $\mathrm{bd}(Y)$ because $\mathbf{x}(t)$ does not approach $\mathrm{bd}(Y)$ as shown in the proof of \textbf{Lemma~\ref{Lemm:ExistLemma1}}: $V$ approaches infinity as $\hat{\mathbf{x}}(t)$ tends to $\mathrm{bd}(Y)$.
Therefore, the limit point belongs to the feasible region.\qed
\end{prf}

From this lemma, we see that there exists the limit point $\hat{\mathbf{x}}^{\infty}$ that satisfies $V(\hat{\mathbf{x}}^{\infty}) = 0$ in the demand- and supply-feasible region.
Moreover, $\hat{\mathbf{x}}^{\infty} = \hat{\mathbf{x}}^{*}$ from the property (ii) in \textbf{Theorem~\ref{Theo:LyapuFunc}}.
This means that there exists an equilibrium state in the feasible region.\qed

\section{Proof of the uniqueness}\label{Sec:App-Unique}
Suppose that multiple equilibria exist.
Consider two different equilibrium states $\hat{\mathbf{x}}^{*}$ and $\hat{\mathbf{y}}^{*}$, which belong to the demand- and supply-feasible region.
Then, substituting these states into [VI-D], we derive the following relationship:
\begin{align}
\sum_{d\in \mathcal{D}}[\mathbf{C}(\mathbf{x}^{*}) - \nabla H_{d}(\mathbf{x}^{*}_{d}) ]\cdot (\mathbf{y}^{*}_{d} - \mathbf{x}^{*}_{d}) \geq 0.\label{Eq:ProofUni_Supposition}
\end{align}
By transforming the left-hand side of this inequality, we have
\begin{align}
&\sum_{d\in \mathcal{D}}[\mathbf{C}(\mathbf{x}^{*}) - \nabla H_{d}(\mathbf{x}^{*}_{d}) ]\cdot (\mathbf{y}^{*}_{d} - \mathbf{x}^{*}_{d}) \\
&= - \left\{\sum_{d\in \mathcal{D}}[ \mathbf{C}(\mathbf{y}^{*}) - \mathbf{C}(\mathbf{x}^{*}) ]\cdot (\mathbf{y}^{*}_{d} - \mathbf{x}^{*}_{d}) 
- \sum_{d\in \mathcal{D}}[ \nabla H_{d}(\mathbf{y}^{*}_{d}) - \nabla H_{d}(\mathbf{x}^{*}_{d}) ]\cdot (\mathbf{y}^{*}_{d} - \mathbf{x}^{*}_{d}) \right\}\notag\\
&\quad - \sum_{d\in \mathcal{D}}[\mathbf{C}(\mathbf{y}^{*}) - \nabla H_{d}(\mathbf{y}^{*}_{d}) ]\cdot (\mathbf{x}^{*}_{d} -  \mathbf{y}^{*}_{d} ).
\end{align}
Owing to the monotonicity of the functions $\mathbf{C}$ and $\nabla H_{d}$, the following relationship holds:
\begin{align}
- \left\{\sum_{d\in \mathcal{D}}[ \mathbf{C}(\mathbf{y}^{*}) - \mathbf{C}(\mathbf{x}^{*}) ]\cdot (\mathbf{y}^{*}_{d} - \mathbf{x}^{*}_{d}) 
- \sum_{d\in \mathcal{D}}[ \nabla H_{d}(\mathbf{y}^{*}_{d}) - \nabla H_{d}(\mathbf{x}^{*}_{d}) ]\cdot (\mathbf{y}^{*}_{d} - \mathbf{x}^{*}_{d}) \right\} < 0.
\end{align}
As $\mathbf{y}^{*}$ is an equilibrium state, we also have
\begin{align}
- \sum_{d\in \mathcal{D}}[\mathbf{C}(\mathbf{y}^{*}) - \nabla H_{d}(\mathbf{y}^{*}_{d}) ]\cdot (\mathbf{x}^{*}_{d} -  \mathbf{y}^{*}_{d} ) \leq 0.
\end{align}
It follows that 
\begin{align}
\sum_{d\in \mathcal{D}}[\mathbf{C}(\mathbf{x}^{*}) - \nabla H_{d}(\mathbf{x}^{*}_{d}) ]\cdot (\mathbf{y}^{*}_{d} - \mathbf{x}^{*}_{d}) < 0.
\end{align}
This contradicts Eq.~\eqref{Eq:ProofUni_Supposition}.
Thus, the supposition must be false and the statement is true.\qed

\section{Efficient solution algorithms}\label{Sec:Algorithm}
As an extended application of the global stability analysis in \textbf{Section~\ref{Sec:ThoreticalAnalysis}}, we present two numerical algorithms for solving the NGEV equilibrium assignment. 
The first algorithm, termed \textit{Algorithm-D}, is based on the NGEV Dynamic-D and the Lyapunov function with disaggregated link flows, whereas the second algorithm, termed \textit{Algorithm-C}, the NGEV Dynamic-C and the Lyapunov function with link costs.
Because these algorithms require repeated evaluation of $\mathbf{F}(\mathbf{c})$, we first explain an efficient (flow-independent) NGEV-based stochastic loading algorithm briefly, referencing \cite{oyama2022markovian}.
We then proceed to provide the detailed description of the Algorithm-D and C.

\subsection{NGEV based dynamic stochastic loading algorithm}
To evaluate the optimal mapping $\mathbf{F}$, the expected minimum cost $\mathbf{s}$ has to be solved (Eq.~\eqref{Eq:OptMap_MinECost}).
As expressed in Eq.~\eqref{Eq:NGEV_ExpMinCost}, the expected minimum cost is formulated in the form of recurrence relation and has the structure of a fixed-point problem: $\mathbf{s}=K(\mathbf{s})$.
Furthermore, the mapping $K$ can be regarded as ``linear" from an abstract algebraic point of view.

To demonstrate this point in detail, we introduce an extended version of the path algebra by \cite{Carre1979-fc}. 
We define the algebra $\mathbb{R}_{\mathrm{path}} \equiv \langle \mathbb{R}_{\epsilon}, \oplus, \otimes \rangle$ as a set $\mathbb{R}_{\epsilon}= \mathbb{R} \cup \{\epsilon\}$ equipped with two binary operations $\oplus$ and $\otimes$.
Letting $[\mathbf{M}]_{ij}$ be the $(i,j)$ element of a matrix $\mathbf{M}$, the operation $\oplus$ for matrices $\mathbf{A}$, $\mathbf{B}\in\mathbb{R}_{\epsilon}^{n\times n}$ is defined in this study as follows:
\begin{align}
    &[\mathbf{A \oplus B}]_{ij} \equiv a_{ij} \oplus_{i} b_{ij},\\
    &\text{where}\quad x\oplus_{i} y 
    \equiv -\cfrac{1}{\theta_{i}} \ln[e^{-\theta_{i}x}+e^{-\theta_{i}y}].
\end{align}
The operation $\otimes$ is also defined as 
\begin{align}
    &[\mathbf{A \otimes B}]_{ij} \equiv (a_{i1} \otimes_{i} b_{1j}) \oplus_{i} \cdots
    \oplus_{i} (a_{ik} \otimes_{i} b_{kj})\oplus_{i}\cdots
    \oplus_{i}(a_{in} \otimes_{i} b_{nj}),\\
    &\text{where}\quad x\otimes_{i} y 
    \equiv -\cfrac{1}{\theta_{i}} \ln[e^{-\theta_{i}x}\cdot e^{-\theta_{i}y}].
\end{align}

By using the algebra $\mathbb{R}_{\mathrm{path}}$, the expected minimum cost can be expressed as follows:
\begin{align}
	\label{eq:mu_path}
	\mathbf{s}_{d} = \mathbf{W} \otimes \mathbf{s}_{d} \oplus \mathbf{e}_{d}	
\end{align}
where $\mathbf{W}\equiv[w_{ij}]\in\mathbb{R}_{\epsilon}^{n\times n}$ is the link weight matrix, and its element is defined as follows:
\begin{align}
	w_{ij}= 
	\begin{cases}
		c_{ij}-\cfrac{1}{\theta_{i}}\log \alpha_{ji},\quad & ij \in L\\
		\infty \quad & ij \not\in L
	\end{cases}\notag
\end{align}
where $\mathbf{e}_{d}$ is the column vector corresponding to node $d$ of the unit matrix $\mathbf{E}\in\mathbb{R}_{\epsilon}^{n\times n}$.

We observe a certain similarity between (\ref{eq:mu_path}) and a system of linear equations $\mathbf{s} = \mathbf{W} \mathbf{s} + \mathbf{e}$ in ordinary matrix algebra. 
It follows that Eq.(\ref{eq:mu_path}) can be solved in an analogous manner to that for the ordinary system of linear equations: if the sequence of powers of $\mathbf{W}$ is convergent, by recursively substituting \eqref{eq:mu_path} into its right-hand side, we have the following equation (herein we omit the superscript $d$ for $\mathbf{s}$ and $\mathbf{e}$ as Eq.~\eqref{eq:mu_path} has the same structure for all $d \in \mathcal{D}$):
\begin{align}
	\label{eq:V}
	\mathbf{s} 
	&= \mathbf{e} \oplus \mathbf{W} \otimes \mathbf{s}\\
	&= \mathbf{e} \oplus \mathbf{W} \otimes [\mathbf{e} \oplus \mathbf{W} \otimes \mathbf{s}] \nonumber\\
	&= [\mathbf{E} \oplus \mathbf{W}] \otimes \mathbf{e} \oplus \mathbf W^{\otimes 2} \otimes \mathbf{s}\notag\\
	& = [\mathbf{E} \oplus \mathbf{W} \oplus \mathbf W^{\otimes 2}] \otimes \mathbf{e} \oplus \mathbf W^{\otimes 3} \otimes \mathbf{s} \nonumber\\
	& = \cdots \nonumber \\
	& = [\mathbf{E} \oplus \mathbf{W} \oplus \mathbf W^{\otimes 2} \oplus \mathbf W^{\otimes 3} \oplus \cdots] \otimes \mathbf{e}.
\end{align} 
That is, solving \eqref{eq:mu_path} reduces to calculating the matrix power series of $\mathbf{W}$ under the algebra $\mathbb{R}_{\rm path}$, which yields the (expected) minimum cost $\mathbf{s}$. 
Once $\mathbf{s}$ is obtained, we can compute the link choice probability matrix $\mathbf{P}$, and finally the link flows $\mathbf{x}$.

\subsection{Algorithms based on the convergent dynamics and Lyapunov functions}

\subsubsection{Algorithm-D: Algorithm based on NGEV Dynamic-D}
The first algorithm minimizes the Lyapunov function $V$ in \textbf{Theorem~\ref{Theo:LyapuFunc}} using a discrete-time version of the disaggregated link dynamic.
The dynamic and the corresponding Lyapunov function $V$ imply that 
\begin{align}
	\mathbf{d}_{d}^{(m)}\equiv \mathbf{y}_{d}^{(m)} - \mathbf{x}_{d}^{(m)},
	\quad
	\text{where}\quad 
	\mathbf{y}_{d}^{(m)}\equiv \mathbf{F}_{d}(\mathbf{C}(\mathbf{x})), 
\end{align}
is a feasible descent direction of the Lyapunov function $V$ when the current solution $\hat{\mathbf{x}}^{(m)}$ is feasible.
This enables us to develop an efficient numerical algorithm for which convergence to the solution of the NGEV-based equilibrium assignment is guaranteed. 
In each iteration of the algorithm, the current solution $\hat{\mathbf{x}}^{(m)}$ is updated by

\begin{align}
	\hat{\mathbf{x}}^{(m+1)}:=
	\hat{\mathbf{x}}^{(m)} + \eta^{(m)}(\hat{\mathbf{y}}^{(m)} - \hat{\mathbf{x}}^{(m)}),
\end{align}
where $\eta^{(m)}$ is the step size and should lie within the interval $0\leq \eta_{m} \leq \eta_{\max}$ to maintain the feasibility of the solution.

It should be noted that a naive line search procedure for determining $\eta_{m}$ should not be used because even a single evaluation of the Lyapunov function $V$ requires performing a NGEV-based network loading (flow-independent traffic assignment). 
To reduce the computational cost for determining the optimal step size, we employ a quadratic interpolation technique similar to that used in the algorithm of \cite{Maher1998-kd} for logit-based SUE assignment. 
To implement the interpolation procedure, the following lemma is useful:
\begin{lemm}
	For a given direction $\hat{\mathbf{d}}\in\mathbb{R}^{L\times D}$, define $\hat{\mathbf{x}}^{\eta}$ as $\hat{\mathbf{x}}^{\eta}\equiv \hat{\mathbf{x}}+\eta \hat{\mathbf{d}}$, and let the associated aggregated vectors in $\mathbb{R}^{L}$ be $\mathbf{x}^{\eta}\equiv \sum_{d\in\mathcal{D}}\mathbf{x}^{\eta}_{d}$, $\mathbf{d}\equiv \sum_{d\in\mathcal{D}}\mathbf{d}_{d}$.
	At any point $(\hat{\mathbf{x}}^{\eta})$, the directional derivative $v(\hat{\mathbf{x}}^{\eta})$ of the function $V(\hat{\mathbf{x}})$ along $\hat{\mathbf{d}}$ is given by 
	\begin{align}
		v(\hat{\mathbf{x}}^{\eta})\equiv
		\mathrm{d}V(\hat{\mathbf{x}}^{\eta})/\mathrm{d}\eta
		= \nabla_{\mathbf{x}}V(\hat{\mathbf{x}}^{\eta})\cdot \mathbf{d}+\nabla_{\hat{\mathbf{x}}}V(\hat{\mathbf{x}}^{\eta})\cdot \hat{\mathbf{d}},
	\end{align}
	where the gradient of the Lyapunov function $V(\hat{\mathbf{x}}^{\eta})$ is given by 
	\begin{align}
		\nabla_{\mathbf{x}}V(\hat{\mathbf{x}}^{\eta})
		= \mathbf{C}(\mathbf{x})-
		[\nabla_{\mathbf{x}}\mathbf{C}(\mathbf{x})](\mathbf{F}(\mathbf{C}(\mathbf{x}))
		- \mathbf{x})
		\quad \text{and}
		\quad 
		\nabla_{\hat{\mathbf{x}}}V(\hat{\mathbf{x}}^{\eta})= -\nabla_{\mathbf{x}_{d}}H_{d}(\mathbf{x}_{d}).
	\end{align}
\end{lemm}

The algorithm based on the NGEV dynamic with disaggregated link flows and the Lyapunov function $V$ can be summarized as follows:

\vspace{2mm}
\begin{algorithm}[H]\label{Algo:StochasticSAlgo}
\caption{\textbf{Algorithm-D}}

{
0. \textbf{Initialization}: 
Set $m := 0$, $\mathbf{x}_{d}^{(0)}:=\mathbf{F}_{d}(\mathbf{C}(\mathbf{0}))$ $\forall d\in\mathcal{D}$, $\mathbf{x}^{(0)}:=\sum_{d\in\mathcal{D}}\mathbf{x}_{d}^{(0)}$.

1. \textbf{NGEV-based assignment for direction finding}: 
\begin{enumerate}
	\item Assign the OD-flow $(\mathbf{q}_{d})_{d\in\mathcal{D}}$ by a NGEV-based assignment procedure based on the current link cost pattern $\mathbf{C}(\mathbf{x}^{(m)})$: 
	\begin{math}
		\mathbf{y}_{d}^{(m)}:=\mathbf{F}_{d}(\mathbf{C}(\mathbf{x}^{(m)}))\ 
		\forall d\in\mathcal{D},
		\quad
		\mathbf{y}^{(m)}:=\sum_{d\in\mathcal{D}}\mathbf{y}_{d}^{(m)}.
	\end{math}
	
	\item Set the direction vector:
	\begin{math}
		\mathbf{d}_{d}^{(m)}:=\mathbf{y}_{d}^{(m)} - \mathbf{x}_{d}^{(m)}\ 
		\forall d\in\mathcal{D},\quad 
		\mathbf{d}^{(m)}:= \mathbf{y}^{(m)} - \mathbf{x}^{(m)}
	\end{math}
	
	\item Set the maximum step-size:
	\begin{align*}
		&[\mathbf{d}_{+}^{(m)}]:=[\mathbf{d}^{(m)}]_{+},\quad 
		\mathcal{L}_{+}:= \{(ij)\in\mathcal{L} \mid d_{ij}^{(m)}>0 \}\\
		&\eta_{max}:=\arg. \max_{ij\in\mathcal{L}_{+}}.
		\left\{
			[\boldsymbol{\mu} - \mathbf{x}^{(m)}]_{ij} / [\mathbf{d}_{+}^{(m)}]_{ij}
		\right\}
		-
		\epsilon.
	\end{align*}
	
	\item Set the directional derivative:
	\begin{align*}
		v_{0}&\equiv v(\hat{\mathbf{x}}^{(m)}):=[\mathbf{C}(\mathbf{x}^{(m)}) - \nabla_{\mathbf{x}}\mathbf{C}(\mathbf{x}^{(m)})\mathbf{d}^{(m)}]\cdot \mathbf{d}^{(m)}
		- \sum_{d\in\mathcal{D}}\nabla_{\mathbf{x}_{d}}H_{d}(\mathbf{x}_{d}^{(m)})\cdot \mathbf{d}_{d}^{(m)}.
	\end{align*}
\end{enumerate}
}

2. \textbf{NGEV-based assignment for curve-fitting}: 
\begin{enumerate}
	\item Set the auxiliary flow for $\eta_{max}$: 
	\begin{math}
		\overline{\mathbf{y}}_{d}^{(m)}:=\mathbf{x}_{d}^{(m)} + \eta_{max}\mathbf{d}_{d}^{(m)},\ 
		\forall d\in\mathcal{D},
		\quad 
		\overline{\mathbf{y}}^{(m)}:=
		\sum_{d\in\mathcal{D}}\overline{\mathbf{y}}_{d}^{(m)}.
	\end{math}
	
	\item Assign the OD-flow $(\mathbf{q}_{d})_{d\in\mathcal{D}}$ by a NGEV-based assignment procedure based on the link cost pattern $\mathbf{C}(\overline{\mathbf{y}}^{(m)})$: 
	\begin{math}
		\mathbf{z}^{(m)}:=\mathbf{F}(\mathbf{C}(\overline{\mathbf{y}}^{(m)}))
	\end{math}
	
	\item Set the auxiliary direction vector: $\mathbf{e}^{(m)}:=\mathbf{z}^{(m)}-\mathbf{x}^{(m)}$.
	
	\item Set the directional derivative:
	\begin{align*}
		v_{1}\equiv v(\overline{\mathbf{y}}^{(m)}):=[\mathbf{C}(\overline{\mathbf{y}}^{(m)}) - \nabla_{\mathbf{x}}\mathbf{C}(\overline{\mathbf{y}}^{(m)})\mathbf{e}^{(m)}]\cdot \mathbf{d}^{(m)}
		- \sum_{d\in\mathcal{D}}\nabla_{\mathbf{x}_{d}}H_{d}(\overline{\mathbf{y}}^{(m)})\cdot \mathbf{d}_{d}^{(m)}.
	\end{align*}
\end{enumerate}

3. \textbf{Determine optimal step size and update the current solution}: 
\begin{enumerate}
	\item Apply quadratic interpolation to estimate the optimal step size:
	\begin{align*}
		\eta^{(m)}:=\min [-v_{0}/(-v_{0}+v_{1}), \eta_{max}].
	\end{align*}
	
	\item Update the current solution: $\mathbf{x}_{d}^{(m+1)}:=\mathbf{x}_{d}^{(m)} + \eta^{(m)}\mathbf{d}_{d}^{(m)}$, $\forall d\in\mathcal{D}$.
\end{enumerate}

4. \textbf{Convergence test}: 
If convergence criterion holds, stop.
Otherwise, set $m:=m+1$ and return to step 1.
\end{algorithm}
\vspace{2mm}

\subsubsection{Algorithm-C: an algorithm based on NGEV Dynamic-C}
The second algorithm minimizes the Lyapunov function $V_{C}$ in \textbf{Theorem~\ref{Defi:NGEVDynamicsC}} by a discrete-time version of the link cost dynamic NGEV Dynamic-C.
That is, the current solution $\mathbf{c}^{(m)}$ is updated by
\begin{align}
	\mathbf{c}^{(m+1)}:=\mathbf{c}^{(m)}+\eta^{(m)}\mathbf{d}^{(m)},
\end{align}
where the feasible descent direction $\mathbf{d}^{(m)}$ of the Lyapunov function is given by 
\begin{align}
	\mathbf{d}^{(m)}:=\mathbf{F}(\mathbf{c}^{(m)}) - \mathbf{C}^{-1}(\mathbf{c}^{(m)}),
\end{align}
which corresponds to the right-hand side of the dynamic. 
Unlike the first algorithm (Algorithm-D), the interpolation technique for determining the step size $\eta^{(m)}$ is not available here because the computation of the gradient $\nabla V_{C}(\mathbf{c})$, which needs the Jacobian matrix $\nabla \mathbf{F}(\mathbf{c})$ of the link demand functions, is prohibitively time-consuming for large scale networks. 
To circumvent this difficulty, we apply the self-regulated averaging scheme proposed in \cite{Liu2009-et}, which is a slight modification of the well-known method of successive averages.
The scheme regulates the increment of $\beta^{(m)}\equiv 1/\eta^{(m)}$ according to the information of the Lyapunov function.
The increment of $\beta^{(m)}$ should be greater than 1, if the iterate tends to diverge, or whenever the Lyapunov function becomes larger; otherwise, the increment of $\beta^{(m)}$ should be less than 1. 
Because $\mathbf{d}^{(m)}$ is always a feasible descent direction of the Lyapunov function $V_{C}$ and the sequence of step size $\{ \eta^{(m)} \}$ satisfies
\begin{align}
	\sum_{m}\eta^{(m)} = +\infty \ 
	\text{and}\ 
	\lim_{m\to +\infty}\eta^{(m)} = 0.
\end{align}
Blum's theorem guarantees that the sequence of $\{ \mathbf{c}^{(m)}\}$ generated by the procedure converges to the equilibrium.
Thus, the second algorithm based on the NGEV dynamic with link costs and the Lyapunov function $V_{C}$ can be summarised as follows:

\vspace{2mm}
\begin{algorithm}[H]\label{Algo:StochasticSAlgo_C}
\caption{\textbf{Algorithm-C}}

0. \textbf{Initialization}: 
Set $m := 0$, $\gamma^{+}>1$, $0<\gamma^{-}<1$, $\beta^{(0)}:=1$, 
\begin{align*}
	\mathbf{c}_{0}:=\mathbf{C}(\mathbf{0}), \ 
	\mathbf{c}^{(1)}:=\mathbf{c}_{0}, \ 
	V^{(0)}=\| \mathbf{C}^{-1}(\mathbf{c}^{(1)})\|^{2}.
\end{align*}

1. \textbf{NGEV-based assignment for direction finding}: 
\begin{enumerate}
	\item Assign the OD-flow $(\mathbf{q}_{d})_{d\in\mathcal{D}}$ by a NGEV-based assignment procedure based on the current link cost pattern $\mathbf{c}^{(m)}$: $\mathbf{y}^{(m)}:=\mathbf{F}(\mathbf{c}^{(m)})$.

	\item Set the direction vector: $\mathbf{d}^{(m)}:=\mathbf{y}^{(m)} - \mathbf{C}^{-1}(\mathbf{c}^{(m)})$, $\tilde{\mathbf{d}}^{(m)}:=[\nabla \mathbf{C}^{-1}(\mathbf{c}^{(m)})]^{-1}\mathbf{d}^{(m)}$.

	\item Set the objective function value: $V^{(m)}\equiv V_{C}(\mathbf{c}^{(m)}):=\| \mathbf{d}^{(m)}\|^{2}$.
\end{enumerate}

2. \textbf{Determine step size and update the current solution}: 
\begin{enumerate}
	\item Determine the step size according to the self-regulated averaging scheme:
	\begin{align*}
		&\text{If}\ V^{(m)}\geq V^{(m-1)}\ \text{then}\ 
		\beta^{(m)}:=\beta^{(m-1)}+\gamma^{+},\ 
		\text{otherwise}\ 
		\beta^{(m)}:=\beta^{(m-1)}+\gamma^{-}; \\
		&\eta^{(m)}:=1/\beta^{(m)};
	\end{align*}
	
	\item Update the current solution: $\mathbf{c}^{(m+1)}:=\mathbf{c}_{0}+[\mathbf{c}^{(m)}+\eta^{(m)}\tilde{\mathbf{d}}^{(m)}-\mathbf{c}_{0}]$.
\end{enumerate}

3. \textbf{Convergence test}: 
If convergence criterion holds, stop.
Otherwise, set $m:=m+1$ and return to step 1.
\end{algorithm}
\vspace{2mm}

\bibliography{NGEV_DTA}
\bibliographystyle{elsarticle-harv}

\end{document}